\documentclass[12pt]{article}
\usepackage{amsmath}
\usepackage{amssymb}
\usepackage{amsfonts}
\usepackage{eucal}
\usepackage[usenames]{color}
\usepackage{graphicx}
\numberwithin{equation}{section}
\newtheorem{Theorem}{Theorem}[section]
\newtheorem{Proposition}[Theorem]{Proposition}

\newtheorem{Remark}[Theorem]{Remark}

\begin{document}
%\title{On the first-passage area of Wiener process with stochastic resetting.}
\title{The first-passage area of Wiener process with stochastic resetting.}
 \author{Mario Abundo\thanks{Dipartimento di Matematica, Universit\`a
``Tor Vergata'', via della Ricerca Scientifica, I-00133 Roma,
Italy. E-mail: \texttt{abundo@mat.uniroma2.it}}  }

\date{}
\maketitle

\begin{abstract}
\noindent
For a one-dimensional Wiener process with stochastic resetting ${\cal X}(t)$, obtained from an underlying Wiener process $X(t),$
we study the statistical properties of  its first-passage time through zero, when starting from $x>0,$ and its first-passage area, that is
the random
area enclosed between the time axis and the path of the process ${\cal X} (t)$ up
to the first-passage time through zero.
By making use of solutions of certain associated ODEs,
we are able to find explicit expressions for the Laplace transforms of the first-passage time and the first-passage area, and their single and joint moments.
\end{abstract}

\noindent {\bf Keywords:} First-passage time, First-passage area, Wiener process, Ornstein-Uhlenbeck process.\\ \\
{\bf Mathematics Subject Classification:} 60J60, 60H05, 60H10.

\section{Introduction }
This paper regards the first-passage area (FPA) of a diffusion process with stochastic resetting; it integrates some other articles
 \cite{abundo:mcap13}, \cite{abvesc17}, \cite{abfu:MCAP19} and  \cite{abundo:SAA21}, in which we studied  the FPA of jump-diffusions, drifted Brownian motion, L\`evy process, and   Ornstein-Uhlenbeck process.
Actually, here we consider a one-dimensional diffusion process in the presence of stochastic resetting, obtained from a underlying diffusion $X(t);$ this kind of process, that we call
Reset Diffusion (RD) process, was first introduced in
\cite{evans2011}, and considered afterwards in \cite{benari}, \cite{denholl}, \cite{evans}, \cite{pal}, \cite{pinsky2023}, \cite{pinsky2020},
%(see also \cite{prebra95})
in particular the corresponding FPA was studied in \cite{singhpal22}, in the case when the underlying process is Brownian motion.
Our aim is to study the statistical properties of the first-passage time (FPT) through zero of a RD process ${\cal X} (t),$ starting from $x>0,$ and its FPA,
namely the area enclosed between the time axis and the path of the process ${\cal X} (t)$ up
to the FPT through zero.
Precisely, in some cases we obtain  explicitly the Laplace transform of FPT and FPA, and their single and joint moments; moreover, we study their behaviors, as
$x \rightarrow 0^+$ and $x \rightarrow + \infty .$ Furthermore,
we provide the
distribution of the maximum displacement of  ${\cal X} (t).$  \par
Really, in the case when the underlying diffusion $X(t)$ is Brownian motion without drift, the FPA was already studied in \cite{singhpal22}, though
the results found therein were obtained making use of special functions.
In contrast, in the present article we utilize nothing but  elementary functions: our arguments are  based on classical results for one-dimensional
diffusions; in fact, the study of the distributions of FPT and FPA is carried out via solutions of certain associated ODEs, also using some results of \cite{singhpal22}. We focus  on the case when the underlying diffusion $X(t)$ is a Wiener process, that is, a Brownian motion with or without drift, but the results can be extended to other RD processes.
\par\noindent
%Really, general formulae for distributions of integral functionals of
%a diffusion with jumps, contained e.g. in [Stochastic Processes by Andrei Borodin], could certainly be applied to diffusion process with stochastic resetting; however, we
%report explicit calculations for our particular case, since we think this can be useful to community.
\par
The study of FPT and FPA of a RD process is important in many areas, e.g. in Biology, in the ambit of stochastic models for the activity of a
neuron subject to resetting (see e.g. \cite{norisa:85} and the references contained therein).
\noindent
FPT and FPA have also important applications in
percolation models and queueing theory (see \cite{dharrama89}, \cite{kea04}, \cite{kearneyetal:05}, \cite{majkea07},  \cite{prebra95});
as an example, in queuing theory the first hitting time to zero is nothing but the
busy period, namely the first instant at which the queue is empty,
and the FPA is the total waiting time  of the ``customers''
during a busy period (see \cite{kea04}). \par\noindent
Other interesting applications of FPA can be found in the ambit of solar activity, dynamics concerning DNA breathing
(see the references in \cite{keamart:14}), and in Economy  (see  \cite{abundo:mcap13});
for a review about FPT and FPA in the case of Brownian motion with resetting, see e.g \cite{maj07}.
\par
Now, we describe precisely the RD process ${\cal X}(t).$ \par\noindent
We consider a one-dimensional temporally homogeneous diffusion process $X(t)$
driven by the SDE:
\begin{equation} \label{diffusion}
dX(t)=  \mu(X(t)) dt + \sigma (X(t)) dB_t ,
\end{equation}
and starting from the position $X(0)=x >0,$
where the drift $\mu (\cdot)$ and diffusion coefficient $\sigma (\cdot)$ are regular enough functions, such that there exists a unique strong solution of the
SDE \eqref{diffusion} for a given starting point,  and $B_t$ is standard Brownian motion; we also assume that the FPT of the diffusion $X(t)$ below zero is finite with probability one. \par\noindent
We suppose that resetting events can occur according to a homogeneous Poisson process with rate $r>0.$
Until the first resetting event $ {\cal R}$ the process ${\cal X}(t)$ coincides with $X(t)$ and it evolves according to \eqref{diffusion} with ${\cal X}(0)=X(0)=x >0;$ when the reset occurs,
${\cal X}(t)$ is set instantly to a position $x_R >0.$ After that, ${\cal X}(t)$ evolves again according to \eqref{diffusion} starting afresh (independently of the past history) from $x_R,$ until the next resetting event occurs, and so on. The inter-resetting times turn out to be independent and exponentially distributed random variables with parameter $r.$ In other words, in any time interval $(t, t+ \Delta t),$ with $\Delta t \rightarrow 0 ^+, $ the process can pass from ${\cal X}(t)$  to the position $x_R$ with probability $r \Delta t  + o( \Delta t),$ or it can continue its evolution according to \eqref{diffusion} with probability $1- r \Delta t + o( \Delta t ).$   The process ${\cal X}(t)$ so obtained is called RD; for any $C^2$ function $f(x),$ its infinitesimal generator is given by (see e.g. \cite{abundo:mcap13}):
\begin{equation} \label{generator}
{\cal L}f(x) = \frac 1 2 \sigma ^2(x) f''(x) + \mu (x) f'(x) +r (f(x_R) -f(x)) \equiv L f(x) +r (f(x_R) -f(x)) ,
\end{equation}
where $Lf(x)= \frac 1 2 \sigma ^2(x) f''(x) + \mu (x) f'(x)$ is the  ``diffusion part'' of the generator, i.e. that regarding the diffusion process $X(t).$
\par
For an initial position $x>0,$ we are concerned with the first-passage time (FPT) of ${\cal X} (t)$ through zero, namely:
\begin{equation}
\tau (x) = \inf \{t>0: {\cal X} (t) = 0 \ | \ {\cal X} (0)=X(0) =x  \} ,
\end{equation}
and the corresponding first-passage area (FPA)
\begin{equation}
A(x) = \int _0 ^ { \tau (x)} {\cal X} (t) dt ,
\end{equation}
which is the area enclosed between the time axis and the path of the process ${\cal X} (t)$ up
to the FPT through zero.
We suppose that both $\tau(x)$ and $A(x)$ are finite with probability one, for any  $x>0.$ Note that
\begin{equation}
\tau (x)=
\begin{cases}
\tau_X (x), \ \ \ \ \ \ \ {\rm if} \ \tau_X (x) < \sigma \\
\sigma + \tau (x_R), \ {\rm if} \ \tau_X (x) \ge \sigma,
\end{cases}
\end{equation}
where $\tau_X (x)$ is the first-hitting time to zero of $X(t)$  starting from $x>0$
and $\sigma$ is an exponentially distributed  random variable with parameter $r>0.$ \par\noindent
Actually, we will limit ourselves to study the case when the underlying process $X(t)$ is a Wiener process, that is Brownian motion with or without drift.
\par
The main qualitative difference between the FPT of the process ${\cal X} (t)$ and the FPT of the underlying diffusion $X(t)$ is that, for the process ${\cal X} (t)$ the moments of the FPT are finite, while for the second one they may be infinite; this is e.g. the case of  Brownian motion starting from $x>0,$ for which as well-known the first-hitting time to zero is finite with probability one, but it has infinite expectation (see e.g. \cite{klebaner}).
\par
The paper is organized, as follows. Section 2 contains some general results, in Section 3 we deal with the case when ${\cal X} (t)$ is Brownian motion with resetting; we find explicit expressions for the Laplace transform and single and joint moments of FPT and FPA;
 moreover, we study the distribution of the maximum displacement of ${\cal X} (t).$ Section 4 contains the analogous results, when $X(t)$ is   Brownian motion with drift $\mu <0.$  In Section 5, we report conclusions and final remarks.

\section{General results}
For $\lambda >0, $ let us consider the Laplace transform (LT) of $\int _0 ^ {\tau (x)} U({\cal X} (s)) ds,$ being  $U(x)=ax+b \ (a, b \ge 0)$ a polynomial of degree one,  that is,
\begin{equation}
M_ \lambda (x) = E \left [ e^ {- \lambda \int _0 ^ {\tau (x)} U({\cal X} (s)) ds } \right ].
\end{equation}
Taking $U(x)=1,$ one obtains the LT of the FPT $\tau (x),$ while for $U(x)=x$ one gets the LT of the FPA $A(x).$ \par\noindent
First, we show that $M_ \lambda (x),$ as a function of $x >0,$  solves a differential equation with boundary conditions. Actually, we can use the analogous result
in \cite{abundo:mcap13}, holding for a general jump-diffusion process, in the special case when the jump governing function is $\gamma (x, u) = x_R - x.$ Thus, we can state the following:
\begin{Proposition} \label{propODEM}
The LT $M_ \lambda (x)$ of $\int _0 ^ {\tau (x)} U({\cal X} (s)) ds$ satisfies the differential problem:
\begin{equation}
\begin{cases}
{\cal L} M_ \lambda (x) = \lambda  U(x) M_ \lambda (x), \\
M_ \lambda (0) =1, \\
\lim _{x \rightarrow + \infty} M_ \lambda (x) = M_ \lambda ( + \infty) < \infty,
\end{cases}
\end{equation}
 or
\begin{equation} \label{ODEforM}
\begin{cases}
LM_ \lambda (x) - \lambda U(x) M_ \lambda (x) - r M _ \lambda (x) + r M_ \lambda (x_R) =0, \\
M_ \lambda (0) =1, \\
\lim _{x \rightarrow + \infty} M_ \lambda (x) = M_ \lambda ( + \infty) < \infty,
\end{cases}
\end{equation}
where $L$ and ${\cal L}$ denote the infinitesimal generator of the underlying diffusion $X(t),$ and of the corresponding process with resetting ${\cal X}(t),$ respectively;
recall  that (see \eqref{generator}) $L$ is defined, for any $C^2$ function $f$  by
\begin{equation}
Lf(x)= \frac 1 2 \sigma ^2 (x) f'' (x) + \mu (x) f'(x),
\end{equation}
and $f'$ and $f''$ denote first and second derivatives of $f.$ \hfill  $\Box$
\end{Proposition}
%\hfill  $\Box$
\begin{Remark}  \label{rem2.2}
As done for  the analogous case in \cite{abundo:mcap13}, the proof of Proposition \ref {propODEM}, can be directly obtained, by using
the following approximation argument; for $\Delta t \rightarrow 0^+,$ we split the integral $\int _0 ^ {\tau (x)} U({\cal X} (s)) ds$ into two addends, by writing:
\begin{equation} \label{formulaM}
M_ \lambda (x) = E \left [ e^ {- \lambda \int _0 ^ {\Delta t} U({\cal X} (s)) ds } e^ {- \lambda \int _{\Delta t} ^ {\tau (x)} U({\cal X} (s)) ds } \right ].
\end{equation}
By setting $u= s- \Delta t$ and
${\cal Y} (u) = {\cal X} (u + \Delta t), \ {\cal Y} (0) = {\cal X} (\Delta t) = x + \Delta x,$
the second integral becomes
$$ \int _ 0 ^ {\tau_1(x)} U( {\cal Y }(u)) du,$$
where $\tau_1(x) = \tau (x) - \Delta t = \tau (x + \Delta x),$ namely it turns out to be equal to
$$ \int _0 ^ { \tau (x+ \Delta x )} U({\cal Y} (u)) du;$$
thus
\begin{equation} \label{iiintegral}
E \left [ e^ {- \lambda \int _{\Delta t} ^ {\tau (x)} U({\cal X} (s)) ds } \right ] = M_\lambda (x+ \Delta x),
\end{equation}
where the expectation in the l.h.s. of \eqref{iiintegral} is calculated under the condition that ${\cal X} (\Delta t) = x + \Delta x .$
For $ \Delta t \rightarrow 0 ^ +,$ at the first order in $\Delta t$ the first integral is equal to $U(x) \Delta t + o (\Delta t),$ and so the first exponential in \eqref{formulaM} is
$1- \lambda U(x) \Delta t + o (\Delta t ).$ \par\noindent
Therefore, at the first order in $\Delta t,$ one obtains
\begin{equation} \label{equation}
 M_ \lambda (x) = (1 - \lambda U(x) \Delta t ) E [ M_ \lambda ( x+ \Delta x)].
 \end{equation}
Now, by conditioning on whether the resetting event $ {\cal R}$ occurs or not in the interval $[t, t+ \Delta t],$  we get:
$$ E[M_ \lambda (x+ \Delta x)] = E[M_ \lambda (x+ \Delta x) \ | {\cal R} ] P({\cal R}) +
E[M_ \lambda (x+ \Delta x) \ | {\cal R} ^C  ] P( {\cal R} ^C).$$
Hence, by using It${\rm \hat{o}}$'s formula, one has:
$$ E[M_ \lambda (x+ \Delta x)] = M_ \lambda (x_R) r \Delta t +  (M_ \lambda(x) + L M_ \lambda (x) \Delta t ) (1-r \Delta t) .$$
Next, by substituting in \eqref{equation} and equating the terms having the same order in $\Delta t,$ we finally obtain
$$LM_ \lambda (x) - \lambda U(x) M_ \lambda (x) - r M _ \lambda (x) + r M_ \lambda (x_R) =0,$$
that is the first of \eqref{ODEforM}. The first boundary condition is due to the fact that, if the process ${\cal X}(t)$ starts from $x=0$, one has $\tau (x)=0;$ as for
the second boundary condition, in absence of resetting one has $ \tau ( + \infty) = + \infty,$ and so $M _\lambda ( + \infty )=0,$ instead, in the presence of resetting, after an exponential time with mean $1/r$ the process is reset to the position $x_R,$ and starting afresh from there, it reaches zero in a finite time ($\tau (x_R)$ is finite with probability one), and so  $M_ \lambda ( + \infty ) < \infty.$
\end{Remark}
\begin{Remark}
Proposition \ref{propODEM} was already proved in \cite{singhpal22} in the case when $X(t)$ is Brownian motion. Note that, for $r=0$ (that is, when no resetting is allowed) one obtains Eq. (2.12) of
\cite{abundo:mcap13}, provided that the jump part in the infinitesimal generator, there denoted by $L_j,$ is set to zero, while the second boundary condition is $M_\lambda ( + \infty )=0.$
\end{Remark}
For any $C^2$ function $f(x),$ we denote by $G$ the operator defined  by:
\begin{equation} \label{definitionG}
G(f(x)) = Lf(x) - rf(x).
\end{equation}
Then, the Eq. in \eqref{ODEforM} can be also written as:
\begin{equation} \label{ODE1forM}
GM_ \lambda (x) = \lambda U(x) M_ \lambda (x) - r M_ \lambda (x_R).
\end{equation}
\noindent If the $n-$th order moment of $\int _0 ^ {\tau(x)}
U({\cal X}(s))ds$ exists  finite, it is provided by:
\begin{equation} \label{generalmoments}
T_ n (x) := E \left [  \left ( \int _0 ^ {\tau(x)}
U(X(s))ds \right ) ^n \right ]  = (-1)^n \left [ \frac {\partial
^n } {\partial \lambda ^n } M_ { \lambda } (x) \right ] _ { \lambda =0} \ , n=1, 2 , \dots .
\end{equation}
By calculating  the $n-$th derivative with respect to $\lambda ,$ at $\lambda =0,$ of
both members  of   \eqref{ODE1forM}, we obtains that, setting $T_0 (x)=1,$
%\begin{Proposition} \label{propositionmoments}
the $n-$th order moments $T_n(x)$  satisfy the ODEs:
\begin{equation} \label{eqmoments}
GT_n (x) = -n U(x) T_{n-1} (x) - r T_n(x_R) , \ x >0,
\end{equation}
with  the constraint  $T_n(0)=0$ and the addition of  an  appropriate  further condition (indeed, we need two conditions to obtain the unique solution of \eqref{eqmoments}).
%\end{Proposition}
%\par \hfill  $\Box$
Note that for $r=0$, \eqref{eqmoments} becomes Eq. (2.19) of \cite{abundo:mcap13};
in particular, for $U(x)\equiv 1,$
 \eqref{eqmoments} is nothing but the celebrated Darling and Siegert's equation  (\cite{darling:ams53}) for the moments of the first-passage time of a diffusion without resetting.
 \par
As regards the joint moments of $\tau(x)$ and $A(x),$ we consider
the joint LT of $\tau(x)$ and $A(x),$ that is $E[e^{- \lambda _1 \tau (x) - \lambda _2 A(x) }] \ ( \lambda _i >0),$ that can be written as
\begin{equation} \label{jointLT}
M_ { \lambda _1, \lambda _2} (x) = E \left [ \exp \left ( - \int _0 ^ { \tau (x)} (\lambda _1 + \lambda _2 {\cal X}(t) ) dt \right ) \right ] .
\end{equation}
As easily seen, one gets:
\begin{equation} \label{momentsfromjointLT}
\frac {\partial M_ { \lambda _1, \lambda _2} (x) }  {\partial \lambda _1 } \mid _{\lambda _1 = \lambda _2 =0}   = - E[ \tau(x)], \
 \frac {\partial M_ { \lambda _1, \lambda _2} (x) }  {\partial \lambda _2 } |_{\lambda _1 = \lambda _2 =0}   = - E[ A(x)],
\end{equation}
and
\begin{equation} \label{jointmomentsfromjointLT}
\frac {\partial ^2 M_ { \lambda _1, \lambda _2} (x)}  {\partial \lambda _1 \partial \lambda _2 } |_{\lambda _1 = \lambda _2 =0}   = E[\tau (x) A(x)].
\end{equation}
By reasoning as done before, taking $U(x)= \lambda _1 + \lambda _2 x,$ we obtain that $M_ { \lambda _1, \lambda _2} (x)$ solves the problem
\begin{equation} \label{problemforjointLT}
\begin{cases}
G M_ {\lambda_1, \lambda _2}(x) =
 (\lambda _1 + \lambda _2 x ) M_ {\lambda_1, \lambda _2}(x) -r M_ {\lambda_1, \lambda _2}(x_R) \\
 M_ {\lambda_1, \lambda _2}(0) =1 \\
 \lim _ { x \rightarrow + \infty }  M_ {\lambda_1, \lambda _2}(x) = M_ {\lambda_1, \lambda _2}(+ \infty) < \infty.
\end{cases}
\end{equation}
Then, by applying $\frac {\partial ^2 } {\partial \lambda _1 \partial \lambda _2 }$ to both members of the first equation in \eqref{problemforjointLT}, and calculating for
$\lambda _1 = \lambda _2 =0,$
we obtain that
$V(x) := E[\tau (x) A(x)]$ is the solution of the problem:
\begin{equation} \label{odeforEtauA}
\begin{cases}
GV(x)= - x E[ \tau(x)] - E[A(x)] -r V(x_R) \\
V(0)=0,
\end{cases}
\end{equation}
with a suitable additional condition. \par\noindent
Note that for $r=0$, \eqref{odeforEtauA} becomes the analogous equation, respectively  obtained in \cite{abvesc17} and \cite{abundo:SAA21}, in the case of drifted Brownian motion and Ornstein-Uhlenbeck process without resetting; of course, now the boundary conditions are different.
\section{Brownian motion with resetting}
In this section, we consider Brownian motion with resetting ${\cal X} (t),$ and we find explicit expressions of the LT and single and joint  moments of FPT and FPA, as solutions of certain differential problems. \par\noindent
Note that for Brownian motion without resetting (i.e. $r=0)$ the moments of the FPT and FPA are infinite (see e.g. \cite{abvesc17}, \cite{kearneyetal:05},
  \cite{maj20}). This follows by the fact that the FPT density of Brownian motion decays as $t^ {-3/2}$ at large time $t.$ In contrast, both the FPT density and the FPA density of  Brownian motion with resetting decay exponentially fast at large values (see e.g. \cite{evans}, \cite{pal}, \cite{singhpal22}), and so the moments of FPT and FPA are finite, in this case.\par\noindent
Since the underlying process of Brownian motion with resetting is $X(t) = B_t,$ then,  for any $C^2$ function $f(x)$ the infinitesimal generator $L$ of $X(t)$ is given  by $Lf(x)= \frac 1 2 f''(x) $
and from \eqref{ODEforM}, we get that the LT of $\int _0 ^ {\tau (x)} U({\cal X} (s)) ds$ solves the problem:
\begin{equation} \label{ODEforMBrownian}
\begin{cases}
\frac 1 2 \frac { \partial ^2 M_ \lambda} {\partial x^2 }(x) - \lambda U(x) M_ \lambda (x) - r M _ \lambda (x) + r M_ \lambda (x_R) =0, \\
M_ \lambda (0) =1, \\
\lim _{x \rightarrow + \infty} M_ \lambda (x) = M_ \lambda ( + \infty) < \infty,
\end{cases}
\end{equation}
Taking $U(x)=1,$ we obtain that the LT of the FPT $\tau(x)$ is the solution of:
\begin{equation} \label{eqforLTtau}
\begin{cases}
\frac 1 2 \frac { \partial ^2 M_ \lambda} {\partial x^2 }(x) - (\lambda +r)  M_ \lambda (x)  + r M_ \lambda (x_R) =0, \\
M_ \lambda (0) =1, \\
\lim _{x \rightarrow + \infty} M_ \lambda (x) = M_ \lambda ( + \infty) < \infty,
\end{cases}
\end{equation}
For fixed $\lambda, r$ and $x_R,$ this is a linear ODE of the second order; the general solution of  the homogeneous equation is
$$ M_ \lambda ^0 (x) = c_1 e^ { - \sqrt {2 (\lambda +r)} \ x } + c_2 e^ { \sqrt {2 (\lambda +r)} \ x },$$
where $c_1, c_2$ are constants with respect to $x$ (they depend only on $r,$ $x_R$ and $\lambda ).$ We search for a particular solution of \eqref{eqforLTtau} in the form
$ \overline M_ \lambda (x) = a = const;$ by imposing this, one finds $a= \frac {r } { \lambda +r} M_\lambda (x_R).$
Then, the general solution of the ODE in the first of \eqref{eqforLTtau} is:
$$ M_ \lambda (x) = c_1 e^ { - \sqrt {2 (\lambda +r)} \ x } + c_2 e^ { \sqrt {2 (\lambda +r)} \ x } + \frac {r } { \lambda +r} M_\lambda (x_R),$$
where
$c_1$ and $c_2$ are constants  to be determined.
From the first boundary condition $M_ \lambda (0) =1,$ it follows
$ c_1 + c_2 + \frac {r } { \lambda +r} M_\lambda (x_R) =1;$ the second boundary condition, i.e. $\lim _{x \rightarrow + \infty} M_ \lambda (x) < \infty $ implies that $c_2 =0,$ and therefore $c_1 = 1 -\frac {r } { \lambda +r} M_\lambda (x_R),$ so we get
\begin{equation} \label{eqparziale}
M _\lambda (x) = E[e^{ - \lambda \tau (x)} ] =\left (1- \frac {r } { \lambda +r} M_\lambda (x_R) \right ) e^ { - x \sqrt {2(\lambda +r)}} + \frac {r } { \lambda +r} M_\lambda (x_R).
\end{equation}
Taking $x=x_R$ in \eqref{eqparziale}, we obtain
\begin{equation} \label{LTtauxR}
M _\lambda (x_R) =  \frac {(\lambda +r)e^{-x_R\sqrt {2(\lambda +r)}} } {\lambda +re^{-x_R \sqrt {2(\lambda +r)}} } ,
\end{equation}
and from \eqref{eqparziale} we finally obtain the explicit expression of the LT of $\tau (x):$
\begin{equation} \label{LTtau}
M_ \lambda (x) = E \left [ e^{- \lambda \tau (x)} \right ]= e^{ -x \sqrt {2(\lambda + r)}} + \frac r {\lambda +r e^{- x_R \sqrt {2(\lambda +r)}}} \left ( e^{- x_R \sqrt {2(\lambda +r)}}-e^{- (x+x_R) \sqrt {2(\lambda +r)}} \right ).
\end{equation}

\begin{Remark} Formula \eqref{LTtau} extends to all $x>0$ Eq. (44) of \cite{singhpal22} that holds for $x=x_R$ (see \eqref{LTtauxR}). \par\noindent
For $r=0,$ namely when no resetting occurs, \eqref{LTtau} provides $E[ e^{- \lambda \tau(x)}] = e^{-x\sqrt{2 \lambda}},$ that is the well-known formula
of the LT of the FPT through zero of BM starting from $x>0;$ this LT corresponds to the inverse Gaussian density for the FPT, i.e.
$f_{\tau (x)} (t)= \frac {x } {\sqrt {2 \pi} } t^{-3/2} e^{- x^2/2t} .$ For $ r \neq 0,$ the LT  $E \left [ e^{- \lambda \tau (x)} \right ]$ given by \eqref{LTtau} can be inverted, at least for $x=x_R$  (see Eq. (46) of \cite{singhpal22} and the comments therein).
\end{Remark}

As regards the LT of $A(x),$ by taking
$U(x)=x$ in \eqref{ODEforMBrownian}, we obtain that the LT of $A(x)$ is the solution of:
\begin{equation} \label{ODEforLTFPABrownian}
\begin{cases}
\frac 1 2 \frac { \partial ^2 M_ \lambda} {\partial x^2 }(x) - r  M_ \lambda (x) = \lambda x M_ \lambda (x)  - r M_ \lambda (x_R), \\
M_ \lambda (0) =1, \\
\lim _{x \rightarrow + \infty} M_ \lambda (x) = M_ \lambda ( + \infty) < \infty,
\end{cases}
\end{equation}
Unfortunately, Eq. \eqref{ODEforLTFPABrownian} cannot be solved in terms of elementary functions; for $x=x_R$ its solution, $M_ \lambda (x_R)= E[e^ {- \lambda A(x_R)}],$ can be written in terms of special functions, precisely the Scorer's  and Airy functions (\cite{abramo}) (see Eqs. (56) and (57) of \cite{singhpal22}). \par\noindent
For $r=0,$ that is for BM without resetting, Eq. \eqref{ODEforLTFPABrownian} is
the Schrodinger equation for a quantum particle
moving in a uniform field (see e.g. \cite{kearneyetal:05}) and it
was solved for any $x >0$ by Kearney and Majumdar (\cite{kearneyetal:05}) in terms of the Airy function (see also Eq. (2.8) of \cite{abvesc17}, and Eq. (2.9) therein, as regards the corresponding density of the FPA).
\subsection{Moments of the FPT}
Now, we go to calculate the first two moments of $\tau (x),$ by solving the ODE \eqref{eqmoments} with $U(x)=1,$ and $T_n(x) = E[ \tau ^n (x)], n=1, 2.$ \par\noindent
As for the mean of $\tau (x),$ taking $n=1,$ we have that $T_1 (x)= E[ \tau (x)]$  is the solution of the problem:
\begin{equation}
\begin{cases}
\frac 1 2 T_1 '' (x) - r T_1 (x) = -1 -rT_1 (x_R) \\
T_1 (0)=0, \  T_1(+ \infty) < + \infty ,
\end{cases}
\end{equation}
Note that the appropriate additional condition is $T_1( \infty) < \infty .$ In fact, since
the FPT density of  Brownian motion with resetting decays exponentially fast at large time (see e.g. \cite{evans}, \cite{pal}, \cite{singhpal22}), for every $x >0$ one has
$E( \tau(x)) < \infty ;$ moreover, from considerations analogue to those of Remark \ref{rem2.2} (final part), it follows that $E(\tau ( + \infty )) < + \infty .$
\par
\begin{figure}
\centering
\includegraphics[height=0.30 \textheight]{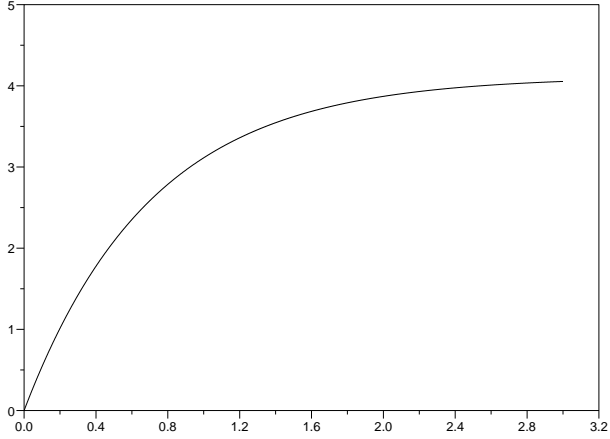}
\caption{Graph of $E[ \tau (x)]$ given by \eqref{Etau}, as a function of $x >0$ for $r=x_R =1$ (on the horizontal axes  $x).$
}
\label{figEtau}
\end{figure}

\noindent
As easily seen, the ODE has general solution
\begin{equation}
T_1(x) = c_1 e^{- x \sqrt {2r}} + c_2 e^{x \sqrt {2r}} + \frac 1 r + T_1 (x_R),
\end{equation}
where $c_1$ and $c_2$ are constants with respect to $x$  to be determined (they depend only on $r$ and $x_R).$
From the first boundary condition $T_1(0) =0,$  we obtain
$ c_1 + c_2 + 1/r + T_1(x_R) =0,$ while
%taking $x=x_R,$ we get
%$$ 0= c_1 e^{- x_R \sqrt {2r}} + c_2 e^{x_R \sqrt {2r}} + \frac 1 r ,$$
from $T_1( + \infty) < + \infty$
we get $c_2 =0.$
Thus, we conclude that $c_1= -1/r - T_1(x_R),$ and so
\begin{equation} \label{Etauparz}
T_1 (x)= \frac { 1 +r T_1 (x_R) } r (1-e^{- x \sqrt {2r}} );
\end{equation}
taking $x=x_R,$ we  obtain:
\begin{equation} \label{EtauxR}
T_1(x_R) = \frac 1 r ( e^{x_R \sqrt {2r}} -1 ),
\end{equation}
and finally from \eqref{Etauparz}:
\begin{equation} \label{Etau}
T_1(x) = E[ \tau (x)] =  \frac 1 r e^ {x_R \sqrt {2r}} \left ( 1- e^{-x \sqrt {2r}} \right ).
\end{equation}
Note that $E[\tau ( + \infty)] = \frac 1 r e^{x_R \sqrt {2r}} < + \infty .$
\begin{Remark}
Of course  formula \eqref{Etau} can be also obtained from  \eqref{LTtau}, by
calculating minus the derivative of $M _ \lambda (x)= E[e^{- \lambda \tau (x)}]$ with respect to $\lambda,$ at $\lambda =0.$
This confirms again that the appropriate additional condition must be $T_1 ( + \infty) < \infty.$
Formula \eqref{Etau} extends to all $x>0$ Eq. (45) of \cite{singhpal22} that provides $E [\tau(x_R)]$ (see \eqref{EtauxR}).
Note that, letting $r$ go to zero in  Eq. \eqref{Etau} it follows $E[\tau (x)] = + \infty,$ which matches the well-known result for BM (see e.g. \cite{abvesc17}).
\end{Remark}
In the Fig. \ref{figEtau}, it is shown an example of the shape of $E[ \tau (x)],$ given by \eqref{Etau}, as a function of $x >0$ for $r=x_R =1.$
\bigskip

As regards the second order moment  of $\tau (x),$ by taking $n=2$ and $U(x)=1$ in \eqref{eqmoments}, we get that $T_2 (x)= E[ \tau ^2 (x)]$  is the solution of the problem:
\begin{equation}
\begin{cases}
\frac 1 2 T_2 '' (x) - r T_2 (x) = - \frac 2 r e^{x_R \sqrt {2r}} (1- e^{-x \sqrt {2r}} ) -r T_2 (x_R) \\
T_2 (0)=0, \  T_2(+ \infty) < + \infty.
\end{cases}
\end{equation}

As easily seen, the ODE has general solution
\begin{equation} \label{generalsolT2}
 T_2(x) = c_1 e^{- x \sqrt {2r}} + c_2 e^{x \sqrt {2r}}- \frac {2x} {r \sqrt {2r}} e^{-(x-x_R) \sqrt {2r}} + \frac 2 { r^2} e^{x_R \sqrt {2r}} + T_2 (x_R) ,
 \end{equation}
where $c_1$ and $c_2$ are constants with respect to $x$ to be determined (they depend only on $r$ and $x_R).$
From the first boundary condition $T_2(0) =0,$  we get
\begin{equation} \label{firstconditionEtau2}
 c_1 + c_2 + \frac 2 {r^2} e^{x_R \sqrt {2r}} + T_2(x_R) =0;
\end{equation}
moreover, from $T_2( + \infty) < + \infty$ we get $c_2=0,$ and so $c_1= -\frac 2 {r^2} e^{x_R \sqrt {2r}} - T_2(x_R).$ \par\noindent
By taking $x=x_R$ in \eqref{generalsolT2} with $c_1= -\frac 2 {r^2} e^{x_R \sqrt {2r}} - T_2(x_R),$ we obtain
\begin{equation} \label{Etau2xRcalculated}
E[\tau ^2(x_R)] = e^{x_R \sqrt {2r}} \left ( \frac 2 {r^2} e^{x_R \sqrt {2r}} - \frac 2 {r^2}  - \frac { 2x_R} { r\sqrt {2r}} \right ) .
\end{equation}
Finally, from \eqref{generalsolT2}:
$$ T_2(x)= E[ \tau ^2(x)] = e^{x_R \sqrt {2r}} \left ( \frac 2 {r^2} e^{x_R \sqrt {2r}} - \frac 2 {r^2}  - \frac { 2x_R} { r\sqrt {2r}} \right ) \left (1- e^{-x \sqrt {2r}} \right )$$
\begin{equation} \label{Etau2x}
-\frac 2 r  e ^{(x_R -x) \sqrt {2r}}  \left ( \frac 1 r + \frac x { \sqrt {2r}}   \right ) + \frac 2 { r^2} e^{x_R \sqrt {2r}}.
\end{equation}
This formula can be also found by calculating from \eqref{LTtau} $\frac {\partial ^2 } {\partial \lambda ^2 } M_\lambda (x),$ at $\lambda =0;$ it confirms that the appropriate additional condition must be $T_2 ( + \infty) < \infty .$
\par
\begin{figure}
\centering
\includegraphics[height=0.30 \textheight]{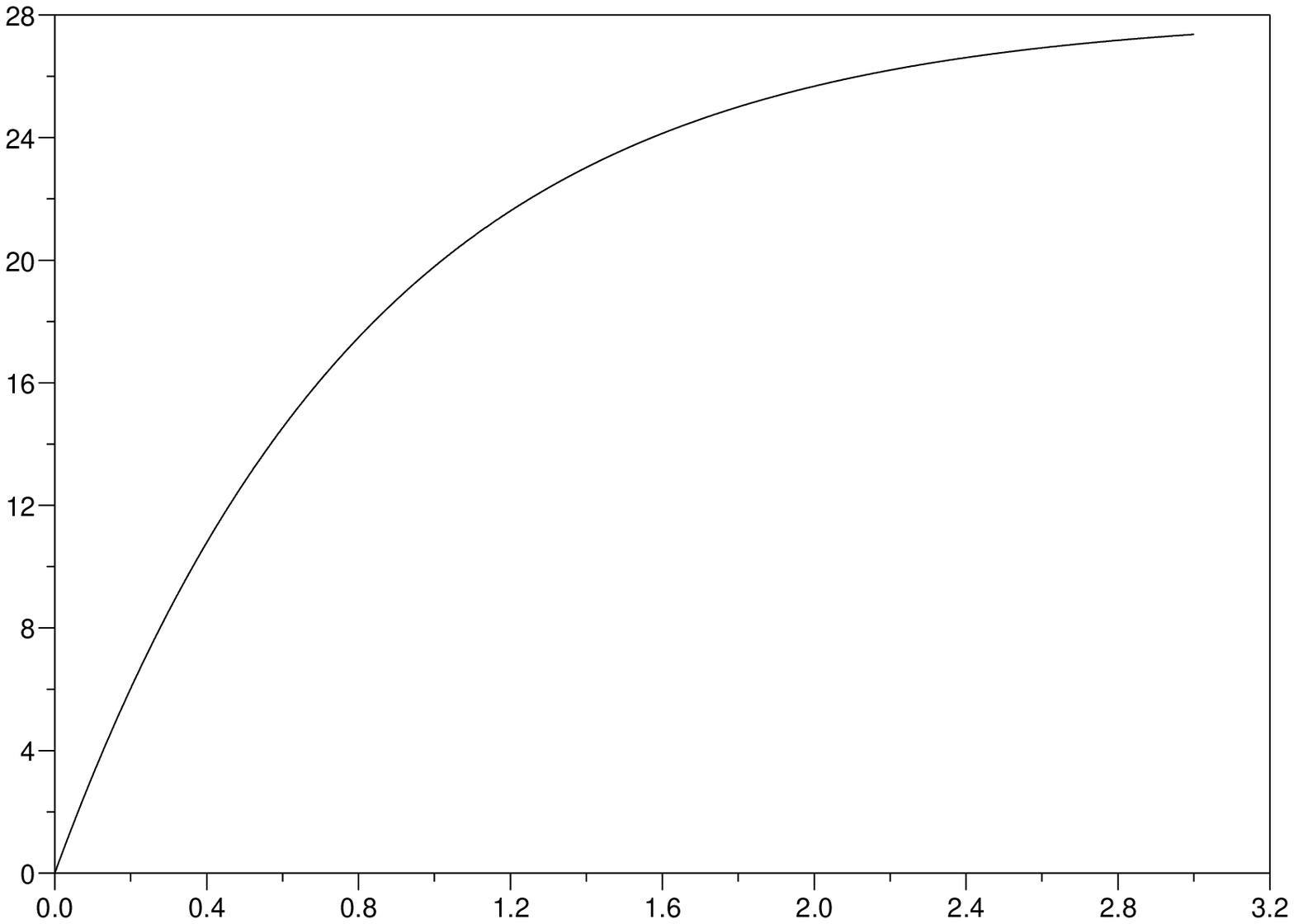}
\caption{Graph of $E[ \tau ^2 (x)]$ as a function of $x >0$ for $r=x_R =1$ (on the horizontal axes  $x).$
}
\label{figEtau2}
\end{figure}

\begin{figure}
\centering
\includegraphics[height=0.30 \textheight]{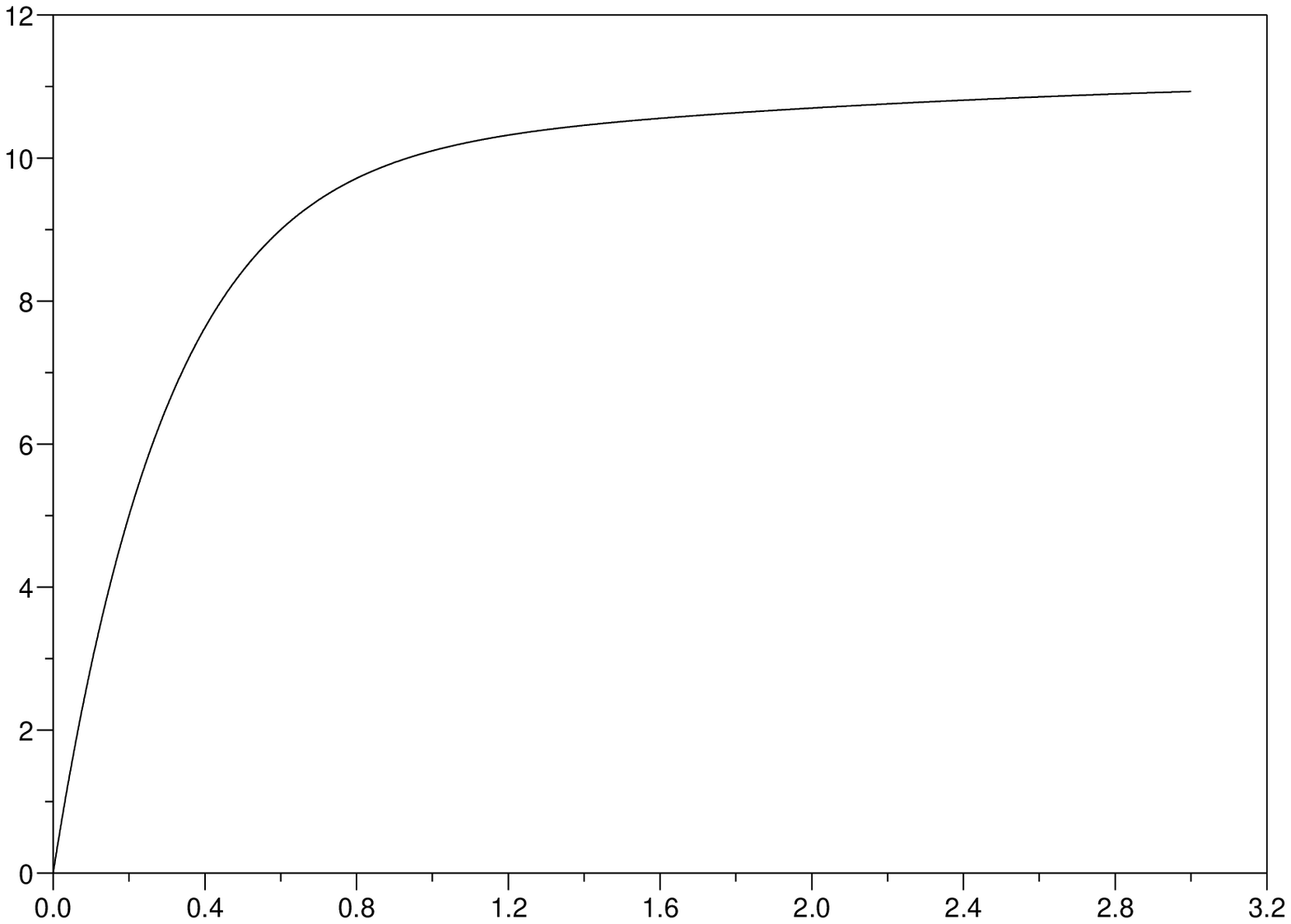}
\caption{Graph of $Var[\tau(x)]$ as a function of $x >0$ for $r=x_R =1$ (on the horizontal axes  $x).$
}
\label{figvartau}
\end{figure}

\noindent
In the Fig. \ref{figEtau2}, it is shown an example of the shape of $E[ \tau ^2(x)],$ given by \eqref{Etau2x}, as a function of $x >0$ for $r=x_R =1.$
In the Fig. \ref{figvartau}, it is shown  the shape of $Var[\tau(x)] = E[ \tau ^2(x)] - E^2[\tau(x)],$  as a function of $x >0$ for $r=x_R =1.$

\par
As regards the behaviors of the first two moments of the FPT $\tau (x)$ for $x \rightarrow 0 ^+$ and $x \rightarrow + \infty,$ from \eqref{Etau} we get
\begin{equation} \label{tauatzero}
E[ \tau (x)] =  a_1 x + o(x), \ {\rm as} \ x \rightarrow 0^+ ,
\end{equation}
where the constant $a_1 = \frac  { \sqrt {2r} } {r } e^{x_R \sqrt {2r}}$ depends on $r$ and $x_R,$
and
\begin{equation} \label{tauatinfty}
\lim _{x \rightarrow + \infty } E [\tau (x)] = a_1 / \sqrt {2r}  .
\end{equation}
From \eqref{Etau2x} we get that
\begin{equation} \label{tausquareatzero}
E[ \tau ^2 (x)] = a_2 x + o(x), \ {\rm as} \ x \rightarrow 0^+,
\end{equation}
where $a_2=  \sqrt {2r} E[\tau ^2 ( x_R)] + \frac { \sqrt 2}{r \sqrt {r}} e^{x_R \sqrt {2r}}$ and $E[\tau ^2 (x_R)]$ is given by \eqref{Etau2xRcalculated};
moreover,
\begin{equation} \label{tausquareatinfty}
\lim _{x \rightarrow + \infty } E[ \tau ^2 (x)] = a_3 .
\end{equation}
where  $a_3 = E[\tau^2 (x_R)] + \frac 2 {r^2} e^{x_R \sqrt {2r}},$
Thus, we obtain
\begin{equation} \label{vartauatzero}
Var [\tau (x)] = a_2  x + o(x), \ {\rm as} \ x \rightarrow 0^+,
\end{equation}
and
\begin{equation} \label{vartauatinfty}
\lim _ {x \rightarrow + \infty} Var [\tau (x)] = a_4 ,
\end{equation}
where $a_4= E[\tau^2 (x_R)] + \frac 2 {r^2} e^{x_R \sqrt {2r}} - \frac 1 {r^2} e^{2x_R \sqrt {2r}}.$

\subsection{Moments of the FPA}
As shown in  \cite{evans}, \cite{maj20}, \cite{pal}, \cite{singhpal22}, the FPA density of Brownian motion with resetting (i.e. for $r>0)$ decays exponentially fast at large values, so the moments of the FPA are finite.
The Laplace transform of the FPA, $E[ \exp(- \lambda A(x)],$ was obtained in \cite{singhpal22} for $x=x_R$ in terms of special functions, and
the first two moments of the FPA at $x=x_R$ were there obtained, by calculating the first and second derivative of $E[ \exp(- \lambda A(x_R)]$  with respect to $\lambda$
at $\lambda =0$ (see Eq. (58) and (59) therein); precisely, it results:
\begin{equation} \label{EAxR}
E[A(x_R)] =  \frac {x_R} r  e^ { x_R \sqrt {2r}} := \alpha _1.
\end{equation}
\begin{equation} \label{truevalueA2xR}
E[A^2(x_R)] = \frac {2 e^ {x_R \sqrt {2r} }} r \left [ \frac {x_R^2} r \left ( \frac 3 4 + e^ {x_R \sqrt {2r} } \right ) + \frac 1 {r^2} - \frac { x_R^3} { 2 \sqrt {2r}}  \right ] -\frac 2 {r^3} := \alpha _2.
\end{equation}
Now, we go to calculate the first two moments of the FPA $A(x)$ for every $x>0,$  by solving the ODE \eqref{eqmoments} with $U(x)=x,$ and $T_n(x) = E[A ^n (x)], n=1, 2 .$ \par\noindent
As regards the mean of $A(x),$ taking $n=1,$ we have that $T_1 (x)= E[A(x)]$  is the solution of the problem:
\begin{equation} \label{ODEformeanA}
\begin{cases}
\frac 1 2 T_1 '' (x) - r T_1 (x) = -x -rT_1 (x_R) \\
T_1 (0)=0, \ T_1(x_R) = \alpha _1,
\end{cases}
\end{equation}
where $\alpha _1$ is given by \eqref{EAxR}; note that $A( + \infty) = + \infty,$ so in contrast with the case of the mean FTP the appropriate additional condition
is $T_1(x_R) = \alpha _1.$
The general solution of the ODE in \eqref{ODEformeanA} is:
\begin{equation} \label{eqforEAparzial}
T_1(x) = c_1 e^ { -x \sqrt {2r}} + c_2 e^ { x \sqrt {2r}} + \frac x r + \alpha _1,
\end{equation}
where $\alpha _1$ is given by \eqref{EAxR}, and $c_1$ and $c_2$ are constants with respect to $x$ to be determined (they depend only on $r$ and $x_R).$ \par\noindent
From $T_1(0)=0,$ we get $c_1 + c_2 + \alpha _1 =0;$ the second condition $T_1(x_R) = \alpha _1$ implies
$c_1 e^{- x_R \sqrt{2r}}+ c_2 e^{ x_R \sqrt{2r}} + \frac {x_R} r =0.$ By solving the system so obtained for $c_1$ and $c_2$
we obtain $c_1 = - \frac 1 r  (e^{ x_R \sqrt{2r}} x_R)$ and    $c_2=0$.
Therefore, finally we get:
\begin{equation} \label{EAx}
T_1(x) =  E[A(x)] = \frac {x_R} r  e^ { x_R \sqrt {2r}} \left ( 1- e^ {- x \sqrt {2r}} \right ) + \frac x r .
\end{equation}
\begin{Remark}
Formula \eqref{EAx} extends to all $x>0$ Eq. (58) of \cite{singhpal22} which provides $E [A(x_R)]$ (see \eqref{EAxR}).
Note that, for $r=0$ from Eq. \eqref{EAx} it follows $E[A(x)] = + \infty,$ which matches the well-known result for BM (see e.g. \cite{abvesc17}).
\end{Remark}
For $x \rightarrow 0^+,$ one has $E[A(x)] =  \left ( \frac {x_R \sqrt {2r}} r e^{x_R \sqrt {2r}} +  \frac 1 r \right )x + o(x),$ while for large positive $x$ it holds
$  E[A(x)] \sim \frac x r.$
\bigskip

\begin{figure}
\centering
\includegraphics[height=0.30 \textheight]{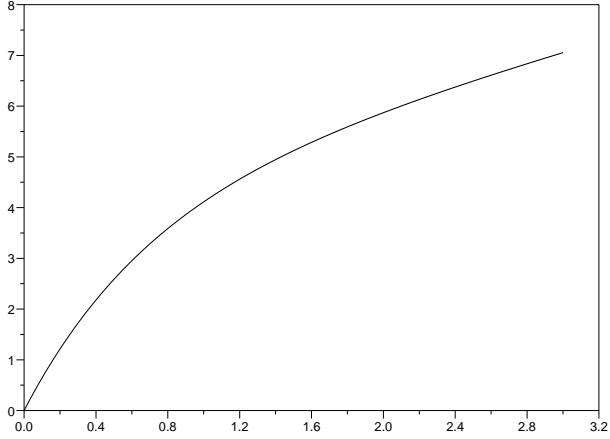}
\caption{Graph of $E[A (x)]$ as a function of $x >0$ for $r=x_R =1$ (on the horizontal axes  $x).$
}
\label{figEAx}
\end{figure}

\noindent
In the Fig. \ref{figEAx}, it is shown an example of the shape of $E[A(x)],$ given by \eqref{EAx}, as a function of $x >0,$ for $r=x_R =1.$

\bigskip

As regards the second order moment of $A(x),$ by taking $n=2$ in
\eqref{eqmoments} with $U(x)=x,$
 we get that $T_2 (x)= E[A^2(x)]$  is the solution of the problem:
\begin{equation} \label{ODEforEA2}
\begin{cases}
\frac 1 2 T_2 '' (x) - r T_2 (x) = -2x E[A(x)] -rT_2 (x_R) \\
T_2 (0)=0, \ T_2(x_R) = \alpha _2,
\end{cases}
\end{equation}
where $\alpha _2$ is given by \eqref{truevalueA2xR} and $E[A(x)]$ is given by \eqref{EAx} (note we have taken the additional condition \eqref{truevalueA2xR}).
\par\noindent
Two independent solutions of the homogeneous equation associated to \eqref{ODEforEA2} are $u_1(x)= e^{- x \sqrt {2 r}}$ and $u_2(x)= e^{ x \sqrt {2r}},$ so its general solution is
$T_2(x) = c_1 u_1(x) + c_2 u_2(x),$
where $c_1$ and $c_2$ are constants with respect to $x$ to be determined (they depend only on $r$ and $x_R).$
By standard methods, we find that
a particular solution of \eqref{ODEforEA2} is
$$ \bar T_2(x) = \frac {2x^2} {r^2} +\frac {2xx_R} { r^2} e^{x_R \sqrt {2r} } + \frac 2 {r^3} + \alpha_2 -
\frac {x_Rx} r e^{-(x-x_R)\sqrt {2r}} \left (\frac x {\sqrt{2r}} + \frac 1 {2r} \right ) .$$
Thus, the general solution of \eqref{ODEforEA2} is:
$$ T_2(x)= c_1 e^{- x \sqrt {2r}} + c_2 e^{ x \sqrt {2r}} + \bar T_2(x) = c_1 e^{- x \sqrt {2r}} + c_2 e^{ x \sqrt {2r}} $$
\begin{equation} \label{generalsolforA2}
 + \frac {2x^2} {r^2} + \frac {2xx_R} {r^2} e^{x_R \sqrt {2r} } +\frac 2 {r^3} + \alpha_2 - e^{-x \sqrt {2r}} \left ( \frac { x_R x^2} { r\sqrt {2r}} e^{x_R \sqrt {2r}} + \frac { x_Rx} { 2 r^2} e^{x_R \sqrt {2r}} \right ) .
\end{equation}
By using the  boundary condition $T_2(0)=0$ and the additional condition  $T_2(x_R) = \alpha _2,$  one gets  a system in the unknown $c_1$ and $c_2,$ whose solution is
$c_1 = - \frac 2 {r^3} - \alpha _2 ,$ and $c_2 =0.$
Then, inserting into \eqref{generalsolforA2}, one finally obtains:
$$
T_2(x) = E[ A^2(x)] = \left ( 1 - e^{-x \sqrt {2r}} \right ) \left ( \alpha_2 + \frac 2 { r^3} \right ) $$
\begin{equation} \label{finaleqforEA2}
+ \frac { 2x} {r^2} \left ( x + x_R e^{x_R \sqrt {2r}} \right ) - \frac {x_R x } r \left ( \frac x {\sqrt {2r}} +
\frac 1 {2r} \right ) e ^{-(x-x_R ) \sqrt {2r}}.
\end{equation}
Eq. \eqref{finaleqforEA2} with $\alpha_2$ given by \eqref{truevalueA2xR} extends to all value of $x >0$ the formula found in \cite{singhpal22} for the expectation of $A^2(x_R)$ (see Eq. (59) therein). \par
\begin{figure}
\centering
\includegraphics[height=0.30 \textheight]{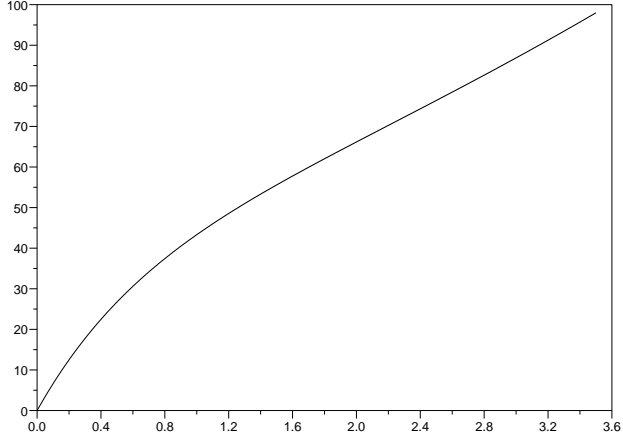}
\caption{Graph of $E[A ^2(x)]$ as a function of $x >0$ for $r=x_R =1$ (on the horizontal axes  $x).$
}
\label{figEA2x}
\end{figure}

\begin{figure}
\centering
\includegraphics[height=0.30 \textheight]{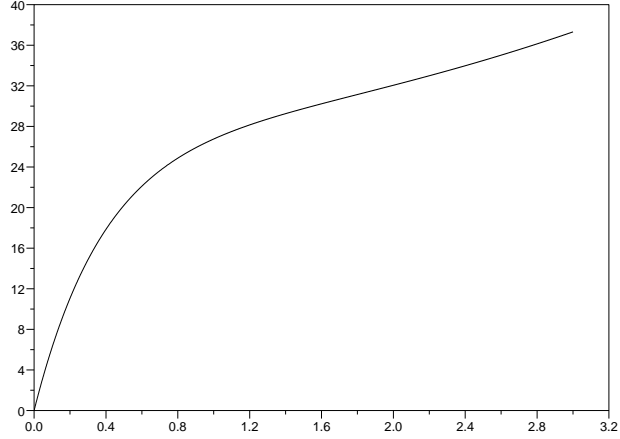}
\caption{Graph of $Var[A(x)],$ as a function of $x >0$ for $r=x_R =1$ (on the horizontal axes  $x).$
}
\label{figvarAx}
\end{figure}

\noindent
In the Fig. \ref{figEA2x}, it is shown the shape of $E[A^2(x)],$ given by \eqref{finaleqforEA2} as a function of $x >0,$ for $r=x_R =1.$ Note that the graph of $E[A^2(x)]$ is not globally concave or convex, but it presents an
inflection point.
In the Fig. \ref{figvarAx}, it is shown the shape of $Var[A(x)],$  as a function of $x >0,$ for $r=x_R =1.$ \par\noindent
From \eqref{finaleqforEA2} it follows that, for $ x \rightarrow 0^+, \ E[A^2(x)] = \left [ \sqrt {2r} \left ( \alpha_2 + \frac 2 {r^3} \right ) +\frac {3 x_R } { 2r^2} e^{x_R \sqrt {2r}}\right ] x +o(x),$
and $Var[A(x)]$ has the same behavior, at the first order in $x;$
moreover, for large $x >0,$
$ E[A^2(x)] \sim \frac 2 {r^2} x^2 ,$
and $  Var[A(x)] \sim   \frac {x^2} { r^2}. $
 \par\noindent

\begin{Remark}
From the previous calculations,
we conclude that the following growth conditions (upper bounds) hold:
\begin{equation} \label{Aatinfty}
E[A(x)] \le const \cdot x , \ {\rm for \ large} \ x >0 ,
\end{equation}
and
\begin{equation} \label{A2atinfty}
E[A^2(x)] \le const ' \cdot x^2 , \ {\rm for \ large} \ x >0.
\end{equation}
Note that similar bounds holds for Brownian motion with negative drift $\mu$ (without resetting, i.e. $r=0),$ because in that case the moments of $A(x)$ grow at most polynomially in $x$ (see \cite{abvesc17}). \par
\end{Remark}

\subsection{Joint moment of $A(x)$ and $\tau(x)$}
In this subsection we find an explicit expression for $E[A(x) \tau(x)],$ i.e. the joint moment of $A(x)$ and $\tau (x).$ The joint LT of $\tau(x), A(x)$ is
$$
M_{\lambda _1, \lambda _2} (x)= E \left [ e^ { - \lambda _1 \tau (x)} e ^ { - \lambda _2 A(x)} \right ], \lambda _1, \lambda _2 >0 .
$$
By reasoning as before (see Proposition 2.1), we find that it satisfies the differential equation:
\begin{equation} \label{eqforjointLT}
(L- r) M_{\lambda _1, \lambda _2} (x) = ( \lambda _1 + \lambda _2 x) M_{\lambda _1, \lambda _2} (x) - r M_{\lambda _1, \lambda _2} (x_R).
\end{equation}
Then, by taking
$\frac {\partial ^2} {\partial \lambda _1 \partial \lambda _2 } $ in both members of Eq. \eqref{eqforjointLT} and calculating it at $\lambda _1 = \lambda _2 =0,$ we obtain that $V(x):= E[\tau(x) A(x)]$ satisfies the differential problem:
\begin{equation} \label{eqforVx}
(L - r) V(x) = -x E[ \tau (x)] -E[A(x)] -r V(x_R) , \ V(0)= 0,
\end{equation}
with a suitable additional condition.
Actually, since $E[A(x) \tau (x)] \le \sqrt {E[A^2(x)]} \sqrt {E[\tau^2(x)]},$ by using \eqref{tausquareatinfty} and \eqref{A2atinfty}, one obtains that the additional condition for the ODE \eqref{eqforVx} must be:
\begin{equation} \label{secondaddcond}
V(x)= E[\tau (x) A(x)] \le const \cdot x, \ {\rm for \ large} \ x >0.
\end{equation}
By standard methods, we find that the general solution of \eqref{eqforVx} is:
$$ V(x)= c_1 e^ { - x \sqrt {2r} } + c_2 e^ {  x \sqrt {2r} } +  \bar V(x), $$
where $c_1$ and $c_2$ are constants with respect to $x$ to be determined (they depend only on $x_R$ and $r),$ and
\begin{equation} \label{barV(x)}
\bar V(x)= \frac {(e^{x_R \sqrt {2r}} +1) x} {r^2} + \frac { x_R e^{x_R \sqrt {2r}}} {r^2}- e^ {-(x-x_R)\sqrt {2r}} \left (\frac {x^2} {2r \sqrt {2r}} +
\frac {(1+2 x_R \sqrt {2r} )x } {4 r^2 }  \right )+ V(x_R)
\end{equation}
is a particular solution of \eqref{eqforVx}.
By imposing the condition $V(0)=0,$ one gets
\begin{equation} \label{firstconditionEAtau}
 0= c_1 + c_2 + \bar V(0)
 \end{equation}
By using \eqref{secondaddcond}, we find that $c_2$ must be zero, so  $c_1 = -\frac {x_R e^{x_R \sqrt {2r}} } {r^2 } - V(x_R),$
and we obtain:
\begin{equation} \label{V(x)parziale}
V(x)= - \left [ \frac {x_R e^{x_R \sqrt {2r}} } {r^2 } + V(x_R) \right ] e^{-x \sqrt {2r}} + \bar V(x).
\end{equation}
By taking $x=x_R$ in \eqref{V(x)parziale}, one obtains an equation in the unknown $V(x_R)$ whose solution provides:
\begin{equation} \label{V(xR)}
V(x_R)= e^{x_R\sqrt {2r}} \left [ \left (8 e^{x_R \sqrt {2r}} -1 \right ) \frac {x_R } {4r^2 } - \frac {3 x_R^2 } {2r \sqrt {2r} } \right ] .
\end{equation}
Finally, the expression for $V(x)= E[\tau(x) A(x)]$ follows from \eqref{V(x)parziale}, \eqref{barV(x)} and \eqref{V(xR)}.
\par\noindent
By using   \eqref{Etau}, \eqref{EAx} and \eqref{V(xR)}, the expressions of  $Cov[A(x), \tau(x)]: = V(x)- E[ \tau(x)] E[A(x)]$
and of the correlation coefficient
\begin{equation}
\rho _{\tau(x), A(x)} := \frac {Cov [\tau(x), A(x)] } { \sqrt {Var[\tau(x)] Var[A(x)]}} ,
\end{equation}
soon follow.
\par\noindent
As easily seen by examination of the various quantities involved, $\rho _{\tau(x), A(x)} $ turns out to be positive for every $x >0,$
 that is, $\tau(x)$ and $A(x)$ are positively correlated; moreover, $\lim _{x \rightarrow 0^+} \rho _{\tau(x), A(x)}= \rho _0 >0 ,$
and $ \lim _{x \rightarrow + \infty} \rho _{\tau(x), A(x)}= \rho _ \infty < \rho _0 .$
The graph of $\rho _{\tau(x), A(x)} $ increases from $\rho _0$ to a maximum value, after that it decreases to $ \rho _ \infty  < \rho _0,$ as $x \rightarrow + \infty .$
This kind of behavior for $\rho _{\tau(x), A(x)} $ was also observed
for drifted BM without reset (see \cite{abvesc17}) and Onrstein-Uhlenbeck process without reset (see \cite{abundo:SAA21}). \par
\begin{figure}
\centering
\includegraphics[height=0.30 \textheight]{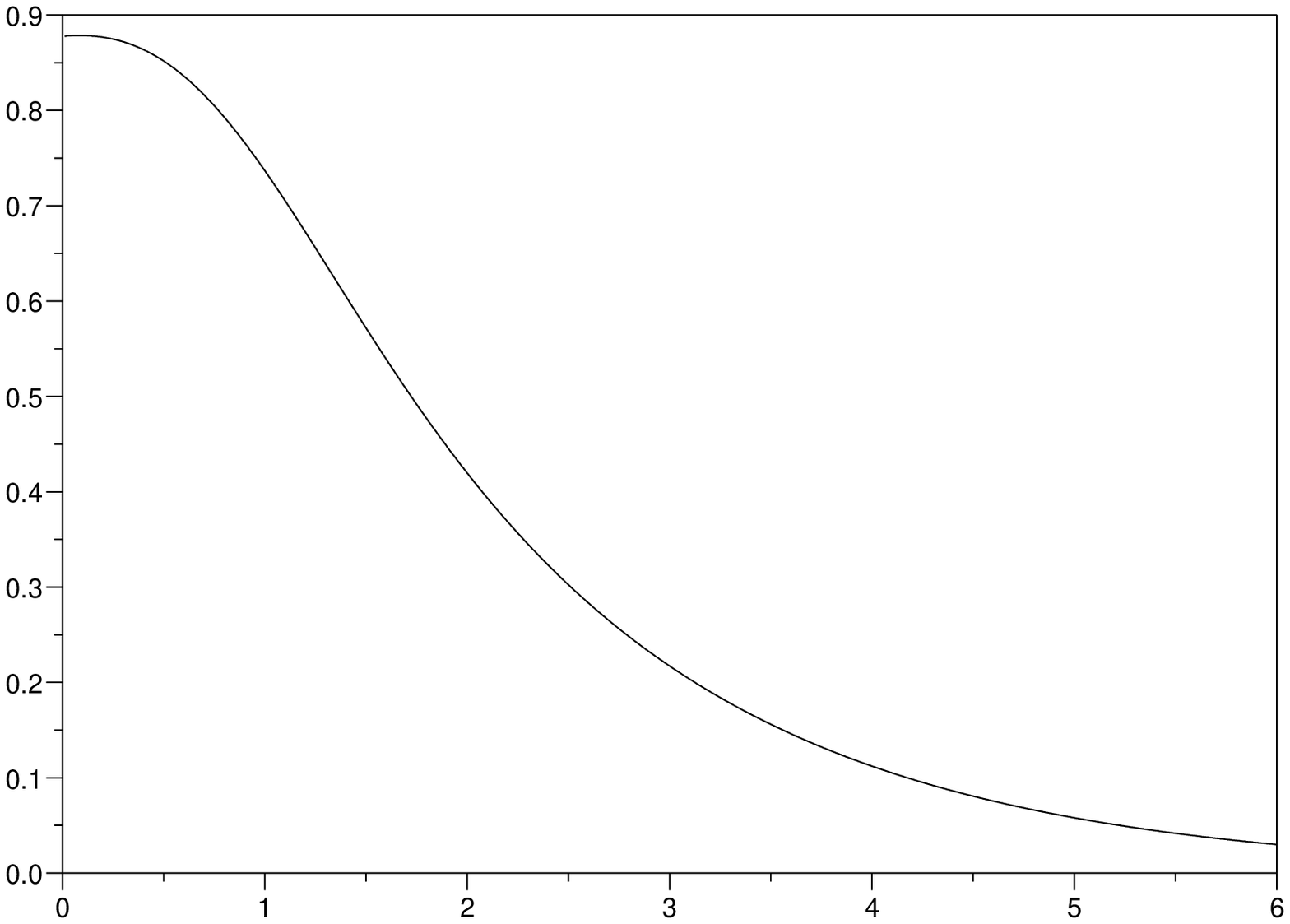}
\includegraphics[height=0.18 \textheight]{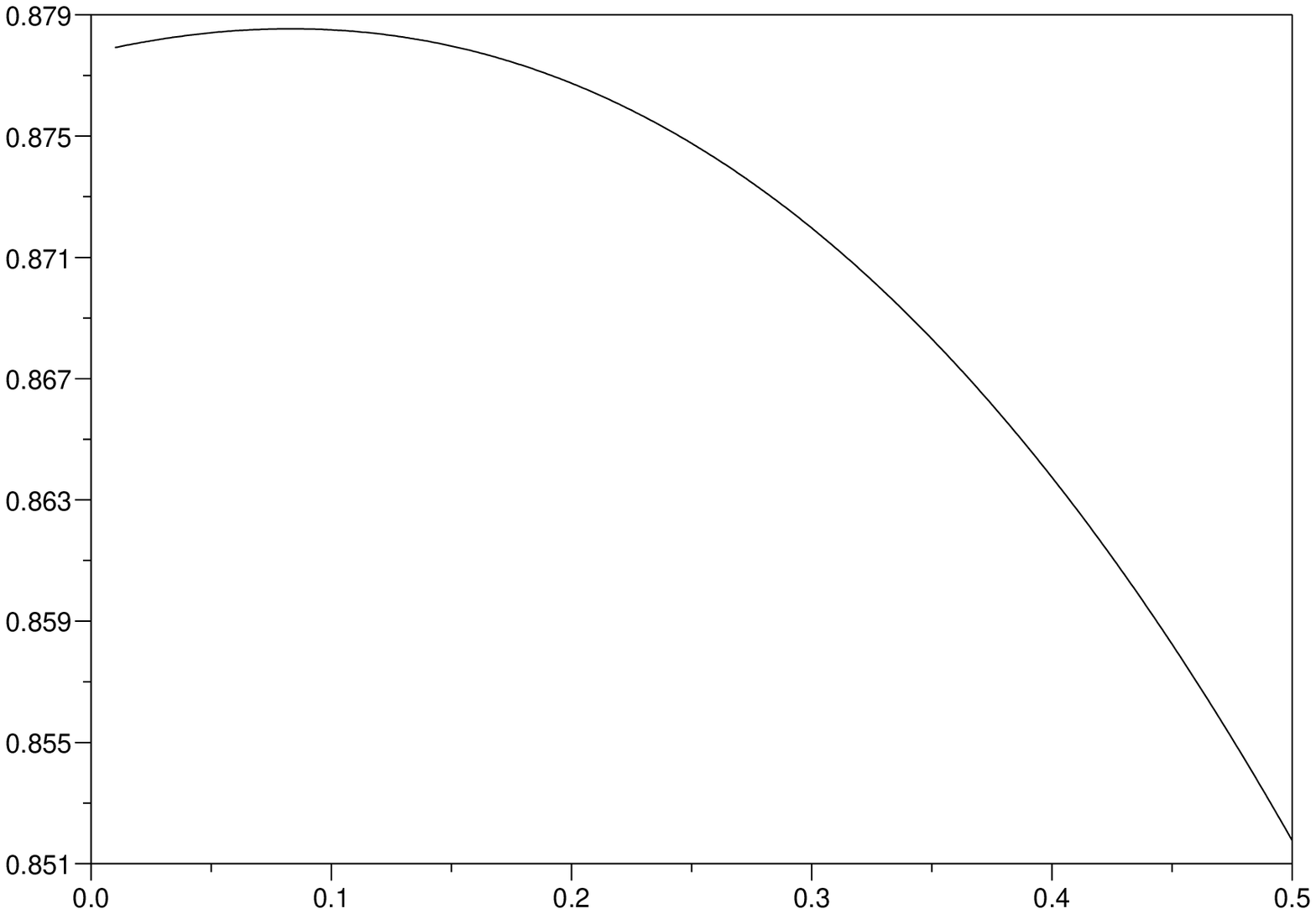}
\caption{Graph of $\rho _{\tau(x), A(x)}$ as a function of $x >0$ for $r=x_R =1$ (on the horizontal axes  $x);$ in the second panel it is shown an enlargement
around $x=0.1 \ .$
}
\label{figrhox}
\end{figure}

\begin{figure}
\centering
\includegraphics[height=0.30 \textheight]{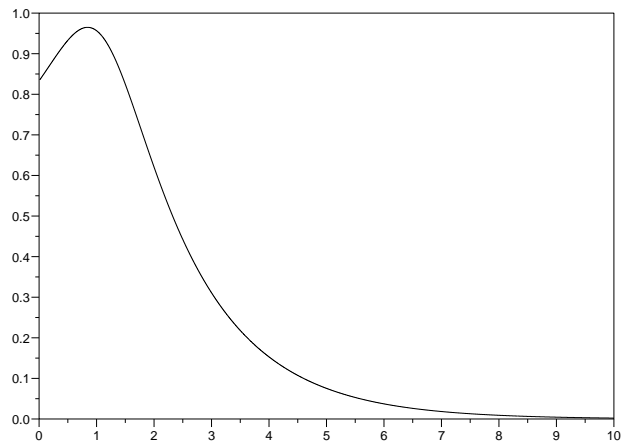}
\caption{Graph of $\rho _{\tau(x), A(x)}$ as a function of $x >0$ for $r=1$ and $x_R =2.$  (on the horizontal axes  $x).$
}
\label{figrhoxbis}
\end{figure}

\noindent
In the Fig. \ref{figrhox}, first panel, we show the graph of $\rho _{\tau(x), A(x)},$ as a function of $x >0$ for $r=x_R =1$ (on the horizontal axes  $x);$ the existence of the maximum is revealed by an enlargement around $x =0.1 \ $ (second panel).
Fig. \ref{figrhoxbis} shows another example of the graph of $\rho _{\tau(x), A(x)},$ obtained for $r=1$ and $x_R =2.$

%\begin{figure}
%\centering
%\includegraphics[height=0.30 \textheight]{grafrho.eps}
%\caption{Graph of $\rho _{\tau(x), A(x)}$ as a function of $x >0$ for $r=x_R =1$ (on the horizontal axes  $x).$
%}
%\label{figrhox}
%\end{figure}

\subsection{Maximum displacement}
\begin{figure}
\centering
\includegraphics[height=0.22 \textheight]{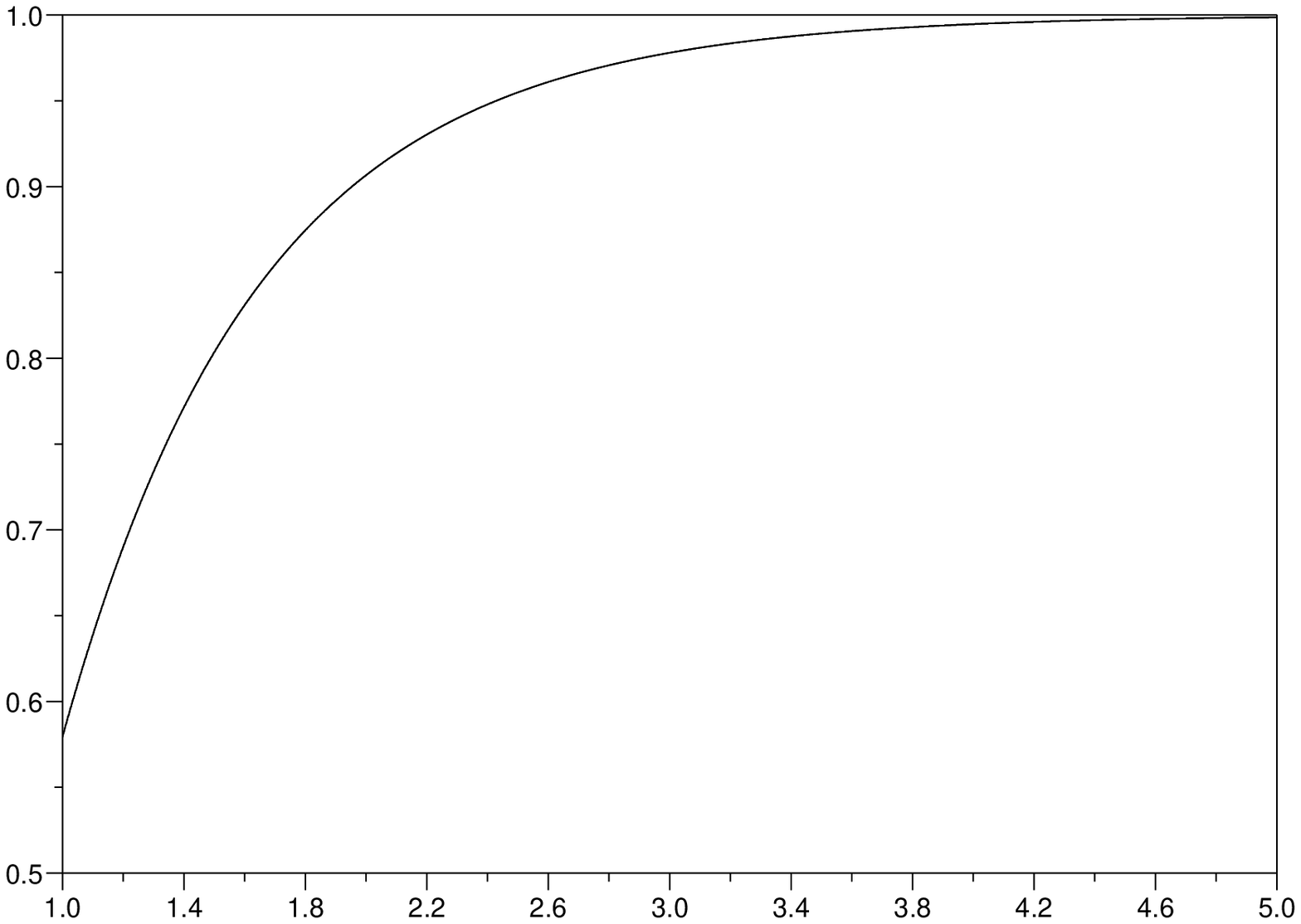}
\includegraphics[height=0.22 \textheight]{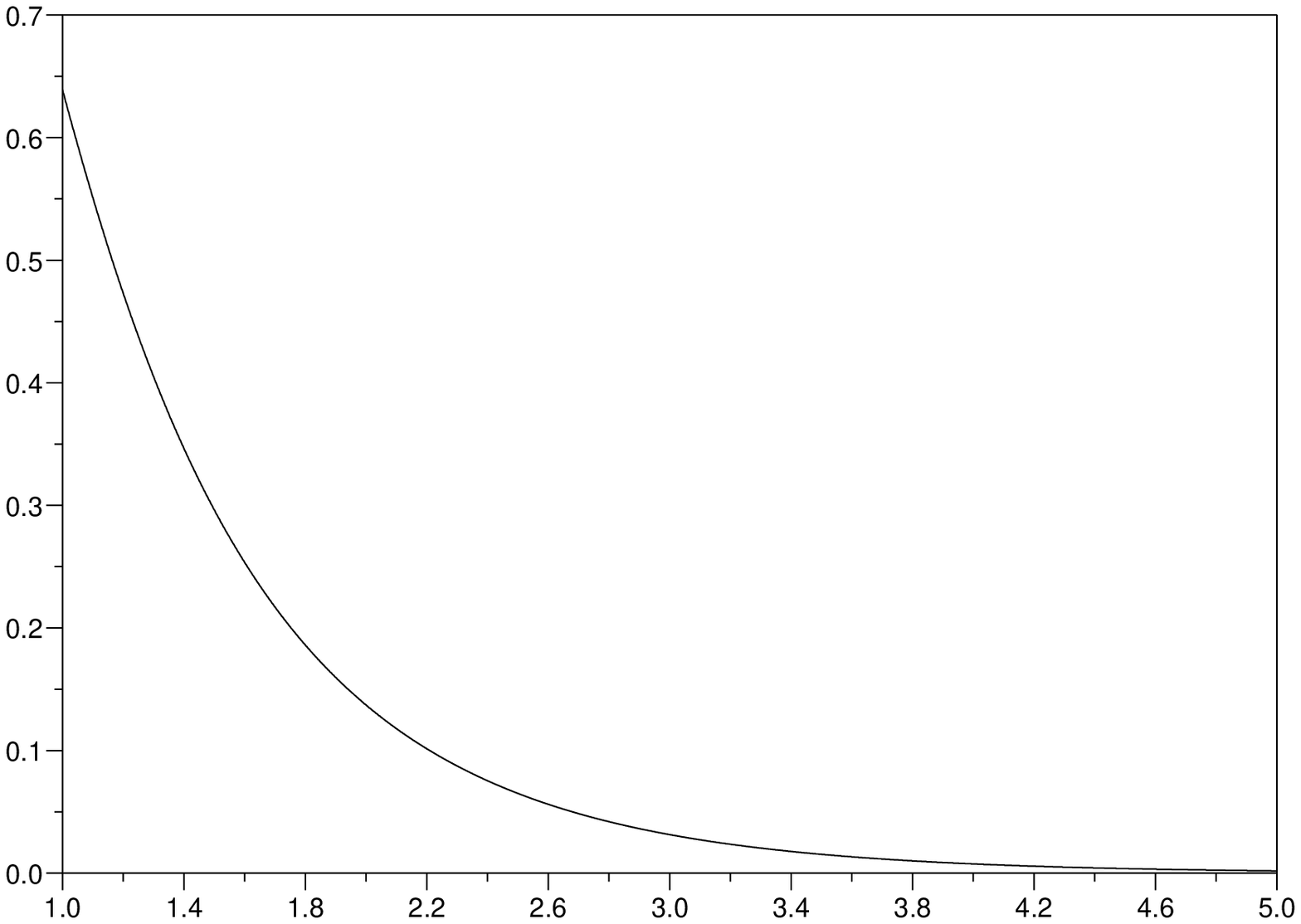}
\caption{Left panel: distribution function of the maximum displacement ${\cal M}_x$ of BM with resetting, for $x=1, x_R=0.5$ and $r=1;$ right panel:
the corresponding density (on the horizontal axes $z \ge x =1).$
}
\label{grafmaxdispl}
\end{figure}

We define the maximum displacement of BM with resetting ${\cal X}(t),$ starting from $x >0,$  as the  r.v.  ${\cal M}_x = \max _ { t \in [0, \tau (x)]} {\cal X}(t)$
(obviously we have ${\cal M}_x \ge x).$
Note that
the event  $\{ {\cal M}_x \le z \}$ occurs if and only if ${\cal X}(t)$ first
exits  the interval $(0, z)$ through the left end $0.$ Therefore, for $z \ge x$  one gets that the distribution function $F_{{\cal M}_x}(z) = P({\cal M}_x \le z )$ is the solution of
the differential equation (see e.g. \cite{abundo:mcap13}, \cite{abundo:SAA21}):
\begin{equation} \label{eqmaxdisplBM}
\begin{cases}
{\cal L} w (x) = \frac { 1 } 2 w''(x) + r w(x_R) - rw(x) =0, \ x\in (0,z) \\
w(0) =1, \ w(z) =0 ,
\end{cases}
\end{equation}
where ${\cal L}$ is the infinitesimal generator of ${\cal X}(t)$ (see  \eqref{generator}).
The general solution of the associated homogeneous equation is $w_0(x)= c_1 e^ { -x \sqrt {2r}} + c_2 e^ { x \sqrt {2r}},$ where for fixed $x$ and $x_R,$ $c_i=c_i(z) $ are indeed functions of $z,$ while a particular solution is $\bar w(x)=a = w(x_R),$ where $a=a(z)$ is a function of $z;$ therefore, the general solution of \eqref{eqmaxdisplBM} is
\begin{equation} \label{partialeqw}
w(x)= c_1  e^ { -x \sqrt {2r}} + c_2 e^ { x \sqrt {2r}} +a,
\end{equation}
where the functions $c_1, \ c_2$ and $a$ are to be found.
By imposing the boundary conditions, and taking $x= x_R$ in \eqref{partialeqw}, we obtain that
$c_1 (z), \ c_2 (z)$ and $a (z)$ must satisfy the system:
\begin{equation}
\begin{cases}
c_1 + c_2 + a =1 \\
c_1 e^ {- z\sqrt{2r}} + c_2 e^ { z\sqrt{2r}} +a =0 \\
c_1 e^ {- x_R \sqrt{2r}} + c_2 e^ { x_R \sqrt{2r}} =0 .
\end{cases}
\end{equation}
By solving the above system, one finds:
\begin{equation} \label{constmaxdisplBM}
\begin{cases}
c_1 = c_1(z) = (e^ {z \sqrt {2r} }-1) ^ {-1} ( e^{ - z \sqrt {2r}} + e^{ - 2x_R \sqrt {2r}}) ^ {-1} \\
c_2= c_2(z)= - e^{ - 2x_R \sqrt {2r}} c_1 \\
a= a(z)= 1- c_1 (z)-c_2 (z).
\end{cases}
\end{equation}
Finally, we get
\begin{equation} \label{maxdispldistr}
F_{{\cal M}_x}(z) = P({\cal M}_x \le z )=
\begin{cases}
0, \ \ \ \ \ \ \ \ \ \ \ \ \ \ \ \ \ \ \  \ \ \ \ \ \ \ \ \ \ \ \ \ \ \ \ z < x \\
c_1 (z) e^ { -x \sqrt {2r}} + c_2 (z) e^ { x \sqrt {2r}} + a(z), \ z \ge x .
\end{cases}
\end{equation}
where $c_i (z)$ and $a(z)= w(x_R)$ are given by
\eqref{constmaxdisplBM}.
As we see,  the distribution function of ${\cal M}_x-x$ appears to have a tail that decays exponentially fast, and so the expectation $E[{\cal M}_x]$ results to be finite.
By calculating the solution of \eqref{eqmaxdisplBM} for $r=0,$ we obtain
that, for BM without resetting the distribution of the maximum displacement is (see also \cite{borodin}):
\begin{equation} \label{maxdispldistrwithoutresetting}
F_{{\cal M}_x}(z) =
\begin{cases}
0, \ \ \ \ \ \ \ z < x \\
1- \frac x z, \ z \ge x ,
\end{cases}
\end{equation}
which implies that the expectation of the maximum displacement of BM without resetting is infinite.
In the Fig \ref{grafmaxdispl} it is shown an example of the distribution function of the maximum displacement ${\cal M}_x$ of BM with resetting, and the corresponding probability density,
obtained by taking the derivative in  \eqref{maxdispldistr},
for $x=1, x_R=0.5$ and $r=1.$

\section{Drifted Brownian motion with resetting}
In this section, we consider drifted Brownian motion with resetting ${\cal X} (t),$ precisely we suppose that the underlying diffusion process is Brownian motion with drift $\mu ,$ that is $X(t)= x + B_t + \mu t,$
and we find explicit expressions of the LT, the moments of FPT and FPA,
and the maximum displacement of  ${\cal X} (t),$ as solutions of differential problems. Since many calculations are analogous to those of the undrifted case, we omit the details, reporting only the results. \par\noindent
Note that the FPT trough zero of  Brownian motion with non-zero drift $\mu $ (without resetting) starting from $x>0$ is finite with probability one, only if the drift  is negative; instead, the
FPT trough zero of drifted Brownian motion with resetting ${\cal X} (t)$ is finite, irrespective of the sign of the drift $\mu$ and  the  moments of FPT are also finite, for any
$0 < r < + \infty$ (see e.g. \cite{pal}).  \par
The infinitesimal generator of $X(t)$ is now given by $Lf(x)= \frac 1 2 \frac { d f^2} {d x^2 } + \mu \frac { d  f} {dx };$ by proceeding as in the case of undrifted Brownian motion with resetting,  the solutions of the various equations cam be  obtained by taking
$\mu + \sqrt{ \mu ^2 + 2(\lambda + r)}$ in place of $\sqrt {2(\lambda +r)}$ in the corresponding formulae.

\subsection{The Laplace transform of $\tau(x)$}
The Laplace transform of $\tau(x)$ turns out to be:
\begin{equation} \label{LTtaudriftedBM}
M_ \lambda (x) = E \left [ e^{- \lambda \tau (x)} \right ]= e^{ -x \left (\mu + \sqrt {\mu ^2 + 2(\lambda + r)} \ \right )} +  M_\lambda (x_R) \frac r {\lambda +r} \left (1- e^{ -x \left (\mu + \sqrt {\mu ^2 + 2(\lambda + r)} \ \right )} \right ),
\end{equation}
where
\begin{equation} \label{LTtauxRdriftedBM}
M _\lambda (x_R) =  \frac {(\lambda +r)e^{-x_R \left (\mu + \sqrt {\mu ^2 + 2(\lambda + r)} \ \right ) }} {\lambda +re^{-x_R \left (\mu + \sqrt {\mu ^2 + 2(\lambda + r)} \ \right ) } }
\end{equation}

\begin{Remark} Formula \eqref{LTtaudriftedBM} appears to be new, because  in  \cite{singhpal22} it was studied only undrifted Brownian motion with resetting.
For $\mu=0$ \eqref{LTtaudriftedBM} and \eqref{LTtauxRdriftedBM} provide again \eqref{LTtau} and \eqref{LTtauxR}.
\end{Remark}

\subsection{The Laplace transform of $A(x)$}
Unfortunately, the corresponding ODE cannot be solved in terms of elementary functions, having non costant coefficients.

\subsection{Moments of the FPT}
The mean of the FPT is:
\begin{equation} \label{Etaudrifted}
E[ \tau (x)] =  \frac 1 r e^ {x_R (\mu +\sqrt {\mu ^2 +2r} ) } \left ( 1- e^{-x (\mu + \sqrt {\mu ^2 + 2r})} \right ).
\end{equation}
For $\mu =0,$ we obtain again \eqref{Etau}. \par
As for the second order moment of  $\tau (x),$  one has:
$$ E[ \tau ^2(x)] $$
$$= e^{x_R (\mu + \sqrt {\mu ^2 + 2r})} \left ( \frac 2 {r^2} e^{x_R (\mu +\sqrt {\mu ^2 +2r} ) } - \frac 2 {r^2}  - \frac { 2x_R} { r\sqrt {\mu ^2 +2r}} \right ) \left ( 1- e^{-x (\mu + \sqrt {\mu ^2 + 2r})} \right )$$
\begin{equation} \label{Etau2xdrifted}
-\frac 2 r  e ^{(x_R -x) ( \mu +\sqrt {\mu ^2 + 2r})}  \left ( \frac 1 r + \frac x { \sqrt {\mu ^2 + 2r}}   \right ) + \frac 2 { r^2} e^{x_R (\mu + \sqrt {\mu ^2 + 2r})}.
\end{equation}
For $\mu =0,$ we obtain again \eqref{Etau2x}.

\subsection{Moments of the FPA}
The mean of the FPA turns out to be:

\begin{equation} \label{EAxdtfted}
 E[A(x)] = E[A(x_R)]  \left ( 1- e^{- x \left (\mu + \sqrt {\mu^2 + 2r} \ \right )} \right ) + \frac x r  + \frac \mu { r^2}
\left ( 1- e^{- x \left (\mu + \sqrt {\mu^2 + 2r} \ \right )} \right )  ,
\end{equation}
where
\begin{equation}
 E[A(x_R)] = e^{ x_R \left (\mu + \sqrt {\mu ^2 + 2r } \ \right )} \left [ \frac { x_R} r + \frac \mu {r^2}  \left (1- e^{- x_R \left (\mu + \sqrt {\mu^2 + 2r} \ \right )} \right ) \right ].
\end{equation}
For $\mu =0,$ we obtain again $E[A(x_R)]= \frac { x_R} r e^{ \sqrt {2r} x_R }$ and $T_1(x)= \frac { x_R} r e^{ \sqrt {2r} x_R } (1-e^{-x \sqrt{2r}}) + \frac x r$
(see \eqref{EAxR} and \eqref{EAx}. \par
As for the second order moment of $A(x),$
in principle, it is possible to find the explicit solution of the corresponding ODE, by proceeding as in the case of undrifted BM with resetting; however the calculations are far heavier and the form of the explicit solution is cumbersome to be written here. Its qualitative behavior can be studied by making use of a software for numerical solutions of linear ODEs of the second order.

\subsection{Maximum displacement}
The distribution function of the maximum displacement ${\cal M}_x$ turns out to be:
\begin{equation}
F_{{\cal M}_x}(z) = P({\cal M}_x \le z )=
\begin{cases} \label{maxdispldistrD}
0, \ \ \ \ \ \ \ \ \ \ \ \ \ \ \ \ \ \ \  \ \ \ \ \ \ \ \ \ \ \ \ \ \ \ \ z < x \\
c_1(z) e^ { d_1x } + c_2 (z) e^ { d_2x} +a (z), \ z \ge x .
\end{cases}
\end{equation}
where
$d_1= - \mu -\sqrt {\mu^2 + 2r} <0, \ d_2 =  - \mu +\sqrt {\mu^2 + 2r} >0,$
and
the functions $c_i (z)$ and $a (z)$ are given by:
\begin{equation} \label{constmaxdisplDBM}
\begin{cases}
c_1 (z)= - \left [e^ {- \left ( \mu + \sqrt {\mu ^2 +2r}\ \right ) z }-1 + e^{-2 x_R \sqrt {\mu ^2 +2r}} \left (1- e^ {- \left ( \mu - \sqrt {\mu ^2 +2r} \ \right ) z } \right ) \right ] ^{-1}
 \\
c_2(z)= - c_1 e^{ - 2x_R \sqrt {\mu ^2 + 2r}} \\
a(z)= 1- c_1(z)-c_2(z).
\end{cases}
\end{equation}
For $\mu =0$ \eqref{maxdispldistrD} becomes \eqref{maxdispldistr}. Once again, the distribution function of
${\cal M}_x-x$ has a tail that decays exponentially fast, and so the expectation $E[{\cal M}_x]$ results to be finite.
Note that for Brownian motion with drift $ \mu <0$ (without resetting) the distribution of the maximum displacement is (see also \cite{borodin}):
\begin{equation} \label{maxdispldistrdriftwithoutresetting}
F_{{\cal M}_x}(z) =
\begin{cases}
0, \ \ \ \ \ \ \ \ \ \ z < x \\
\frac {e^{-2 \mu x} - e^ {-2 \mu z} } {1- e^{-2 \mu z} } , \ z \ge x ,
\end{cases}
\end{equation}
and its density is
\begin{equation} \label{maxdispldensdriftwithoutresetting}
f_{{\cal M}_x}(z) =
\begin{cases}
0, \ \ \ \ \ \ \ \ \ \ z < x \\
\frac { \mu e^{- \mu x} \sinh (\mu x) } {\sinh ^2( \mu z) }, \ z \ge x ,
\end{cases}
\end{equation}
from which it follows that,
unlike the case $\mu =0,$  the expectation of the maximum displacement of BM with drift $\mu <0$  is
%$E({\cal M}_x) = x + \frac {e^{-2 \mu x} -1 } {-2 \mu } \ln \left ( \frac 1 {1- e^{2 \mu x}} \right ) < + \infty$.
$E[{\cal M}_x] = x - \frac {1- e^{-2 \mu x}  } {2 \mu } \ln \left (1- e^{2 \mu x} \right ) < + \infty$.
\par\noindent
Note that, for $\mu \rightarrow 0 ^-$ the last expression tends to $+ \infty$ (see the comment after \eqref{maxdispldistrwithoutresetting}).

\section{Conclusions and final remarks}
We have studied  a one-dimensional Wiener process with stochastic resetting ${\cal X}(t)$, obtained from an underlying drifted Brownian motion  $X(t)= x + \mu t + B_t,$
starting from $x >0,$ where $B_t$ is standard Brownian motion, and the drift $\mu$ is non positive.
We have supposed that resetting events  occur according to a homogeneous Poisson process with rate $r>0;$
until the first resetting event the process ${\cal X}(t)=X(t)$ evolves as a Brownian motion with drift with ${\cal X}(0)=X(0)=x >0;$ when the reset occurs,
${\cal X}(t)$ is set instantly to a position $x_R >0.$ After that, ${\cal X}(t)$ evolves again as Brownian motion with drift, starting afresh (independently of the past history) from $x_R,$ until the next resetting event occurs, and so on.
The inter-resetting times turn out to be independent and exponentially distributed random variables with parameter $r.$
%This means that, in any time interval $(t, t+ \Delta t),$ with $\Delta t \rightarrow 0, $ the process can pass from ${\cal X}(t)$  to the position $x_R$ with probability $r \Delta t  %+ o( \Delta t),$ or it can continue its evolution according to Brownian motion with drift with probability $1- r \Delta t + o( \Delta t ).$ \par
We have studied the statistical properties of  the first-passage time (FPT) of ${\cal X}(t)$ through zero, when starting from $x>0,$ and its first-passage area (FPA)
(i.e. the  random area enclosed between the time axis and the path of the process ${\cal X} (t)$ up to the FPT through zero).
By solving  certain associated ODEs, we have obtained explicitly the Laplace transform of the FPT and the FPA of ${\cal X}(t),$ and their single and joint moments;
moreover, we have provided the
distribution of the maximum displacement of  ${\cal X} (t).$
The calculations have been shown in some detail only in the case $\mu =0;$ for $\mu \neq 0$ we have only reported the corresponding formulae.
\par\noindent
Notice that some results, regarding undrifted Brownian motion with stochastic resetting, were already obtained  in  \cite{singhpal22}, by using special functions.
\par
We emphasize that, as regards the FPT of drifted Brownian motion with resetting ${\cal X}(t)$, the main qualitative difference with the case of the drifted Brownian motion $X(t)$
(without resetting) is that the FPT of ${\cal X}(t)$ through zero
is finite and it has finite expectation for every $x >0,$
whilst the mean of the FPT of $X(t)$ is finite only for $\mu <0,$ being  infinite  for $\mu =0,$ though
the hitting time to zero of Brownian motion starting from $x>0$ is finite with probability one.
Actually,
the FPT density of Brownian motion decays as $t^ {-3/2}$ at large time $t.$ In contrast, both the FPT density and the FPA density of  Brownian motion with resetting decay exponentially fast at large values (see e.g. \cite{evans}, \cite{pal}, \cite{singhpal22}), and so the moments of FPT and FPA result to be finite.
\par
In principle, the arguments of this paper can be used to obtain differential equations for the Laplace transforms, the moments of FPT and FPA and the distribution of the maximum displacement also for other diffusion processes with stochastic resetting, such as
the Ornstein-Uhlenbeck process with stochastic resetting (see e.g. \cite{dubey}), namely  when the underling process is the diffusion $X(t) = X_{OU} (t)$ which is driven by the SDE
$dX(t) = - \nu X(t) + \sigma dB_t,$ for positive constants $\nu$ and $\sigma.$
For instance, if $\tau$ denotes again the FPT of $X_{OU}(t)$ through zero, when starting from $x>0,$ then $T(x)= E[\tau (x)]$ turns out to be solution to:
\begin{equation} \label{etauforOUreset}
\begin{cases}
\frac 1 2 T''(x) - \nu x T'(x) - r T(x) = -1 -r T(x_R) \\
T(0)=0, \ T_1( + \infty) < \infty.
\end{cases}
\end{equation}
The additional condition $T_1( + \infty) < \infty $ can be replaced by  $\lim _{\nu \rightarrow + \infty} T(x) =0.$ \par\noindent
The differential equation above is complicated enough, but it can be solved in principle in the following way.
Let $\psi (x) = E[ e^ {- r \tau _{OU} (x)}] $ be the Laplace transform of $\tau _{OU} (x),$ which is explicitly given by Eq. (3.49) of \cite{abundo:SAA21} in terms of a cylinder parabolic function. Then, we search for a solution to \eqref{etauforOUreset} of the form
$T(x)= \psi (x) v(x),$ where $v(x)$ is a function to be determined.  By calculating first and second derivatives and inserting into \eqref{etauforOUreset} one obtains the ODE:
$$\frac 1 2 \psi (x) v'' (x) + v'(x) ( \psi '(x) - \nu x \psi (x)) =c,$$
with $c= -1 -r T(x_R),$ that can be solved by quadratures, in the unknown function $v(x),$ and so $E[ \tau (x)]= T(x)= \psi (x) v(x)$ can be obtained.
Analogous way can be followed in principle to solve the differential equations for the LT of $\tau(x)$ and the mean of $A(x).$ \par\noindent
Note that in \cite{dubey} the explicit expressions of the LT  of the FPT $\tau(x)$ of Ornstein-Uhlenbeck process with stochastic resetting  and  its mean $E[ \tau (x)]$ were explicitly obtained in terms of special functions.


\bigskip




\begin{thebibliography}{}

\bibitem {abramo}
Abramowitz, M. and Stegun I.A. (1965) \newblock {\it Handbook of Mathematical Functions: with formulas, graphs, and mathematical tables.}
\newblock New York: Dover.

\bibitem {abundo:mcap13}
Abundo, M. (2013) \newblock
On the first-passage area of a one-dimensional jump-diffusion process.
\newblock{\it Methodol Comput Appl Probab}  {15} 85--103.
DOI 10.1007/s11009-011-9223-1.


\bibitem {abvesc17}
Abundo, M. and Del Vescovo, D. (2017) \newblock  On the joint distribution of first-passage time and first-passage area of  drifted Brownian motion.
\newblock {\it Methodol Comput Appl Probab} {19} 985--996. DOI 10.1007/s11009-017-9546-7

\bibitem {abfu:MCAP19}
Abundo, M., and Furia, S. (2019) \newblock
Joint Distribution of First-Passage Time and First-Passage Area of Certain  L\`evy Processes.
\newblock{\it Methodol Comput Appl Probab} 21 1283--1302.  https://doi.org/10.1007/s11009-018-9677-5

\bibitem {abundo:SAA21}
Abundo, M.  (2023) \newblock
The first-passage area of Ornstein-Uhlenbeck process revisited.
\newblock{\it Stochastic Analysis and Applications} 41(2) 358--376. DOI: 10.1080/07362994.2021.2018335.


%\bibitem {abpir18}
%Abundo, M. and Pirozzi, E.
%\newblock Integrated Stationary Ornstein-Uhlenbeck Process, and Double Integral Processes.
%\newblock {\it Physica A} {494} (2018) 265--275.


%\bibitem {abupir:2021}
%Abundo, M. and Pirozzi, E. \newblock
%Fractionally Integrated Gauss-Markov processes and applications.
%\newblock{\it Commun Nonlinear Sci Numer Simulat}  {105862} (2021) 1--20.
%https://doi.org/10.1016/j.cnsns.2021.105862


\bibitem {benari}
Ben-Ari, I. (2012) \newblock
Principal eigenvalue for Brownian motion on a
bounded interval with degenerate
instantaneous jumps.
\newblock   {\it Electron. J. Probab.} 17(87) 1--13. DOI: 10.1214/EJP.v17-1791


\bibitem {borodin}
Borodin, A.N., Salminen, S. (1996) \newblock {\it Handbook of
Brownian Motion-Facts and Formulae.}
\newblock Birkhauser Verlag Basel, Basel.

\bibitem {darling:ams53}
Darling, D. A. and Siegert, A.J.F. (1953) \newblock
The first passage problem for a continuous Markov process.
\newblock   {\it Ann. Math. Statistics} {24} 624--639.

%\bibitem {vescthesis}
%Del Vescovo, D. (2017) \newblock {\it On first-passage time and first-passage area of some one-dimensional diffusion processes.}
% \newblock {Unpublished Math Degree Thesis, Tor Vergata University, Rome, Italy}.


\bibitem {denholl}
den Hollander, F., Majumdar, S.N., Meylahn, J.M., and Touchette, H. (2019) \newblock
Properties of additive functionals of Brownian motion with resetting.
\newblock   {\it J. Phys. A: Math. Theor.} 52 175001. DOI 10.1088/1751-8121/ab0efd

\bibitem {dharrama89}
Dhar, D. and Ramaswamy, R. (1989)  \newblock Exactly solved model of self-organized critical phenomena.
 \newblock {\it Physical Review Letters} 63(16), p. 1659.


\bibitem {dubey}
Dubey, A. and Pal, A. (2023) \newblock First-passage functionals for Ornstein Uhlenbeck process with stochastic resetting,
 \newblock {\it arXiv:2304.05226v1}.

\bibitem {evans}
Evans, M. R., Majumdar, S.N., and Schehr, G. (2020)
\newblock Stochastic Resetting and Applications.
\newblock {\it  J. Phys. A: Math. Theor.} 53 193001.

\bibitem {evans2011}
Evans, M. R., Majumdar, S.N. (2011) \newblock Diffusion with Stochastic Resetting.
\newblock {\it  Phys. Rev. Lett.} 106 160601.


\bibitem {kea04}
Kearney, M.J. (2004) \newblock On a random area variable arising in discrete-time queues and compact
directed percolation. \newblock {\it Journal of Physics A: Mathematical and General}  37(35)  8421.

\bibitem {keamart:21}
Kearney, M.J. and Martin, R.J. (2021) \newblock
Statistics of the first passage area functional for an Ornstein-Uhlenbeck process.
\newblock{\it J. Phys. A: Math. Theor.} {54} 055002 1--14.
DOI 10.1088/1751-8121/abd677

\bibitem {kearneyetal:05}
Kearney, M.J.and Majumdar, S.N. (2005) \newblock
On the area under a continuous time Brownian motion till its first-passage
time. \newblock {\it J. Phys. A: Math.} Gen 38 4097--4104.

\bibitem {keamart:14}
Kearney, M.J., Pye A.J. and Martin, R.J. (2014) \newblock
On correlations between certain random variables associated with first
passage Brownian motion.
\newblock{\it J. Phys. A: Math. Theor.} {47}(22) 225002. doi:10.1088/1751-8113/47/22/225002

\bibitem {klebaner}
Klebaner, F.C. (2006) \newblock {\it Introduction to Stochastic Calculus with Applications.}
\newblock Imperial College Press, London.

%\bibitem {kreutz}
%{Kreutz-Delgado, K. } \newblock
%{Mean Time-to-Fire for the Noisy LIF Neuron - A Detailed Derivation of the Siegert Formula.}
%\newblock
%{{\it 	arXiv:1501.04032}} (2015)

\bibitem {maj20}
Majumdar, S.N. and Meerson, B. (2020) \newblock Statistics of first-passage Brownian functionals.
\newblock{ \it J. Stat. Mech.: Theory and Experiment} (2) 023202.

\bibitem {majkea07}
Majumdar, S.N. and Kearney, M.J. (2007) \newblock Inelastic collapse of a ball bouncing on a randomly
vibrating platform.  \newblock {\it Physical Review E} 76(3) p. 031130.

\bibitem{maj07}
Majumdar, S.N. (2007) \newblock Brownian functionals in physics and computer science. \newblock In The Legacy Of
Albert Einstein: A Collection of Essays in Celebration of the Year of Physics 93--129.


\bibitem {norisa:85}
Nobile, A.G., Ricciardi, L.M., and Sacerdote, L.  (1985) \newblock
Exponential trends of Ornstein-Uhlenbeck first-passage-time densities.
\newblock{\it J. Appl. Prob.} {22} 360--369.

\bibitem{pal}
Pal, A., Chatterjee, R., Reuveni, S. and Kundu, A. (2019) \newblock Local time of diffusion with stochastic
resetting. \newblock {\it Journal of Physics A: Mathematical and Theoretical} 52(26) p. 264002.


\bibitem{pinsky2023}
Pinsky, R.G. (2023) \newblock Large time probability of failure in diffusive search with resetting for a random target in
$R^d-$ a functional analytic approach.\newblock {\it Transactions of the American Mathematical Society}, 376(4).

\bibitem{pinsky2020}
Pinsky, R.G. (2020) \newblock Diffusive search with spatially dependent resetting.
\newblock {\it Stoch. Process. Their Appl.} 130(5) 2954--2973.

\bibitem {singhpal22}
Prashant Singh, and Arnab Pal (2022) \newblock
First-passage Brownian functionals with stochastic
resetting.
\newblock {J. Phys. A: Math. Theor.} 55 234001. https://doi.org/10.1088/1751-8121/ac677c

\bibitem {prebra95}
Prellberg, T. and Brak, R.  (1995) \newblock Critical exponents from nonlinear functional equations for
partially directed cluster models. \newblock {\it Journal of statistical physics} 78(3) 701--730.



\bibitem {riccetal:99}
Ricciardi, L.M., Di Crescenzo, A., Giorno, V. and Nobile, A.G. (1999) \newblock
An outline of theoretical and algorithmic approaches to first passage time problems with applications to biological modeling.
\newblock{\it Math. Japonica} {50}(2) 247--322.

%\bibitem {ross}
%Ross, S.M. \newblock
%Introduction to probability models. Tenth Ed. Burlington: Academic Press, Elsevier, 2010.

\bibitem {siegert:51}
Siegert, A.J.F. (1951) \newblock
On the First Passage Time Probability Problem
\newblock{\it Phys. Rev.} {81}(4)  617--623.

\bibitem {thomas:75}
Thomas, M.U. (1975) \newblock
Some mean first-passage time approximations for the Ornstein-Uhlenbeck process.
\newblock{\it J. appl. Prob.} {12} 600--604.


\end{thebibliography}
\end{document}